\documentclass[12pt,a4paper]{article}
\usepackage[utf8]{inputenc}
\usepackage[english]{babel}
\usepackage{amsmath}
\usepackage{amssymb}
\usepackage{amsthm} 
\usepackage{graphicx}
\usepackage{longtable} 
\usepackage{array}
\newtheorem{theorem}{Theorem}
\newtheorem{lemma}{Lemma}
\newtheorem{remark}{Remark}

\newtheorem{corollary}{Corollary}

\usepackage{multirow}

\usepackage[colorlinks=true,linkcolor=blue,citecolor=blue]{hyperref}
\usepackage[title]{appendix}

\usepackage[svgnames]{xcolor}
\def\ep{\varepsilon}

\def\th{\vartheta}

\usepackage{lineno}

\title{$L_2$-small ball asymptotics \\for Gaussian random functions: a survey}
\author{
Alexander Nazarov\footnote{St.Petersburg Dept of Steklov Institute, Fontanka 27, St.Petersburg, 191023, Russia, 
    and St.Petersburg State University, Universitetskii pr. 28, St.Petersburg, 198504, Russia. E-mail: al.il.nazarov@gmail.com. } \ and \setcounter{footnote}{6}
Yulia Petrova\footnote{Pontifícia Universidade Católica do Rio de Janeiro (PUC-Rio), R. Marquês de São Vicente, 124, Gávea, Rio de Janeiro, 22451-040, Brazil, and St. Petersburg State University, Universitetskii pr. 28, St.Petersburg, 198504, Russia. E-mail: yu.pe.petrova@yandex.ru.}}
\date{\today}

\begin{document}
\maketitle

\begin{abstract}
    This article is a survey of the results on asymptotic behavior of small ball probabilities in $L_2$-norm. 
\end{abstract}
\hfill \textit{Dedicated to the memory of M.S. Birman and M.Z. Solomyak,}

\hfill \textit{founders of St. Petersburg school of spectral theory}
\footnotesize
\tableofcontents
\normalsize

\section{Introduction}

The theory of small ball probabilities (also called small deviation probabilities) 
is extensively studied in recent decades (see the surveys by M.A.~Lifshits~\cite{Lif1999}, W.V.~Li and Q.-M.~Shao~\cite{LiShao2001} and V.R.~Fatalov~\cite{Fat2003}; for the extensive up-to-date bibliography see~\cite{LifshitsAirtable}). Given a random vector $X$ in a Banach space $\mathfrak B$, the relation
\begin{align}
\label{eq:exact-small-ball-probability}
\mathbb{P}\left\{\|X\|_{\mathfrak B}\leq\ep\right\}\sim f(\ep), \qquad \ep\to0,
\end{align}
(here and later on $f(\ep)\sim g(\ep)$ as $\ep\to0$ means $\lim\limits_{\ep\to0}f(\ep)/g(\ep)=1$)
is called an \textit{\textbf{exact}} asymptotics of small deviations. 
Typically, this probability is exponentially small, and often a \textit{\textbf{logarithmic}} asymptotics is studied, that is the relation  
$$
\log\big(\mathbb{P}\{\|X\|_{\mathfrak B}\leq\ep\}\big)\sim f(\ep), \qquad \ep\to0.
$$
Theory of small deviations has numerous applications including the complexity problems~\cite{PapWas1990}, the quantization problem~\cite{GraLusPag2003}, \cite{LusPag2004}, the calculation of the metric entropy for functional sets~\cite[Section~3]{Lif1999}, functional data analysis~\cite{Fer2006}, nonparametric Bayesian estimation~\cite{VaarZan2008},  machine learning~\cite{RasmWil2006}
(more applications and references can be found in~\cite[Chapter~7]{LiShao2001}). 

The discussed topic is almost boundless, so in this paper we focus on the most elaborated (and may be the simplest) case, where $X$ is \textit{\textbf{Gaussian}} and the norm is \textit{\textbf{Hilbertian}}. Let $X$ be a Gaussian random vector in the (real, separable)  Hilbert space $\mathfrak H$ (the basic example is a Gaussian random function in $L_2({\cal O})$, where ${\cal O}$ is a domain in $\mathbb R^d$). We always assume that $\mathbb{E} X=0$ and denote by $\mathcal G_X$ corresponding covariance operator (it is a compact non-negative operator~in~$\mathfrak H$). 

Let $(\mu_k)$ be a non-increasing sequence of positive eigenvalues of $\mathcal G_X$, counted with their multiplicities, and denote by $(\varphi_k)$ a complete orthonormal system of the corresponding eigenvectors.\footnote{If $\mathcal G_X$ has non-trivial null space then $\varphi_k$ form a complete system in its orthogonal complement.} It is well known (see, e.g.,~\cite[Chapter 2]{Lif2012-book}) 
that if $\xi_k$ are i.i.d. standard normal random variables, then we have the following distributional equality\footnote{For the special choice $\xi_k=\frac{1}{\sqrt{\mu_k}}\,(X,\varphi_k)_{\mathfrak{H}}$ this equality holds almost surely.}
$$
X\stackrel{d}{=}\sum\limits_{k=1}^\infty \sqrt{\mu_k}\,\varphi_k \,\xi_k,
$$
which is usually called the \textit{\textbf{Karhunen--Lo{\`e}ve expansion}}.\footnote{Up to our knowledge stochastic functions given by similar series were first introduced by D.D.~Kosambi~\cite{Kos1943}. However, K.~Karhunen~\cite{Karh1946}, \cite{Karh1947} and M.~Lo{\`e}ve~\cite{Loe1948-book} were the first who proved the optimality in terms of the total mean square
error resulting from the truncation of the series. An analogous formula for stationary Gaussian processes was introduced by M.~Kac and A.J.F.~Siegert in~\cite{KacSie1947}.
}
By the orthonormality of the system $(\varphi_k)$ this implies 
\begin{align}
    \label{eq:KL}
    \|X\|_{\mathfrak{H}}^2\stackrel{d}{=}\sum\limits_{k=1}^{\infty}\mu_k\xi_k^2.
\end{align}
Therefore the small ball asymptotics for $X$ in the Hilbertian case is completely determined by the eigenvalues $\mu_k$.

Notice that if $\sum\limits_{k=1}^\infty \mu_k<\infty$ then the series in~\eqref{eq:KL} converges almost surely (a.s.), otherwise it diverges a.s. The latter is impossible for $X\in\mathfrak{H}$, so, in what follows we always assume that $\mathcal{G}_X$ is the trace class operator.

\begin{remark}
Formula~\eqref{eq:KL} shows that if the covariance operators of two Gaussian random vectors $X$ and $Y$ have equal spectra (excluding maybe zero eigenvalues), then their Hilbertian norms coincide in distribution.\footnote{In fact, the converse statement is also true.} In this case $X$ and $Y$ are called \textit{\textbf{spectrally equivalent}}. Trivially, such Gaussian random vectors have the same $L_2$-small ball asymptotics.

Concrete examples of spectrally equivalent Gaussian processes are well known, see, e.g., \cite{DonMarYor1991}, \cite[Example~4]{KlepLeBre2002}, \cite{Pyck2003-phd},~\cite[Theorem~4.1]{Naz2009-non-separated},~\cite{Liu2013-addBM}. Some multivariate generalizations can be found in~\cite{DehPecYor2006}, \cite{PecYor2006},~\cite{Deh2006}.

Recently A.I.~Nazarov and Ya.Yu.~Nikitin~\cite{NazNik2022-spect-equiv} (see also~\cite{NazNik2021-online}) developed
a general operator approach to the problem of spectral equivalence and obtained a number of new examples.
\end{remark}

Our paper is organized as follows. In Section~\ref{sec:firstworks} we describe the first separated results on small ball probabilities in Hilbertian norm. In Section~\ref{sec:secondwave} we formulate two crucial results: the Wenbo Li comparison principle and the Dunker--Lifshits--Linde formula. These results formed the base for the systematic attack of the problem. 

Further progress in the field is mainly based on the methods of spectral theory. In Section~\ref{sec:green} we give an overview of results on exact asymptotics. These results are mostly related to the spectral theory of differential operators. Section~\ref{sec:logarithm} is devoted to the logarithmic asymptotics, corresponding results are related to the spectral theory of integral operators. 

In Appendix we provide the history (up to our knowledge) of results on $L_2$-small ball asymptotics for concrete processes.
\medskip

Let us introduce some notation.

For $\mathcal{O}\subset \mathbb{R}^d$, $\|\cdot\|_{2,\mathcal{O}}$ stands for the norm in $L_2(\mathcal{O})$:
$$
\|X\|_{2,\mathcal{O}}^2=\int\limits_{\mathcal{O}}|X(x)|^2\,dx.
$$

The Sobolev space $W_p^m[0,1]$ is the space of functions $u$ having
continuous derivatives up to $(m-1)$-th order when $u^{(m-1)} $ is
absolutely continuous on $[0,1]$ and $u^{(m)}\in L_p[0,1]$.  For $p=2$, it is a Hilbert space.

A measurable function $f(x)$ is called regularly varying at infinity (see, e.g.~\cite[Chapter~1]{Sen1976-book}), if it is of constant sign on $[A, \infty)$, for some $A>0$, and there exists $\alpha\in\mathbb{R}$ such that for arbitrary $\rho>0$ we have
\begin{equation*}
\lim\limits_{x\to\infty}\frac{f(\rho x)}{f(x)}=\rho^\alpha.
\end{equation*}
Such $\alpha$ is called the index of the regularly varying function. Regularly varying function of index $\rho=0$ is called the slowly varying function (SVF). For instance, the functions $\log^\sigma (x)$, $\sigma\in\mathbb{R}$, are slowly varying at infinity.

A function $f(x)$ is called regularly varying at zero, if $f(\frac{1}{x})$ is a regularly varying at infinity. We say that a sequence $(a_k)$ has regular behavior if $a_k=f(k)$, where $f(x)$ is regularly varying at infinity.

\section{First works}
\label{sec:firstworks}

The oldest results about $L_2$-small ball probabilities concern classical Gaussian processes on the interval~$[0,1]$. R.~Cameron and W.~Martin~\cite{CamMar1944} proved the relation for the Wiener process $W(t)$
\begin{align*}
\mathbb{P}\bigl\{\|W\|_{2,[0,1]}\leq\ep\bigr\}\sim \frac{4\ep}{\sqrt{\pi}}\exp\bigl(-\frac18\,\ep^{-2}\bigr),
\end{align*}
while T.~Anderson and D.~Darling~\cite{AndDar1952} established the corresponding asymptotics for the Brownian bridge~$B(t)$
\begin{align}
\label{eq:BB}
\mathbb{P}\bigl\{\|B\|_{2,[0,1]}\leq\ep\bigr\}\sim \frac{\sqrt{8}}{\sqrt{\pi}}\exp\bigl(-\frac18\,\ep^{-2}\bigr).
\end{align}
The latter formula describes the lower tails of famous Cramer--von Mises--Smirnov $\omega^2$-statistic. 

One should notice that a minor difference (rank one process) between Wiener process and Brownian bridge influences the power term in the asymptotics but not the exponential one.

In general case the problem of $L_2$-small ball asymptotics was solved by G.~Sytaya~\cite{Syt1974}, but in an implicit way. The main ingredient was the Laplace transform and the saddle point technique. 
The general formulation of result by Sytaya in terms of operator theory states as follows.
\begin{theorem}[{\cite[Theorem~1]{Syt1974}}]
    \label{thm:sytaya-general}
    As $\ep\to0$, we have
    \begin{multline}
    \label{eq:sytaya-general}
        \mathbb{P}\Bigl\{\|X\|_{\mathfrak{H}}\leq \ep\Bigr\}
        \sim \bigl(4\gamma^2\pi\,\mathrm{Tr}(\mathcal G_X R_\gamma(\mathcal G_X))^2\bigr)^{-\frac 12}\times\\
        \times\exp\biggl(-\int\limits_0^\gamma \mathrm{Tr} (\mathcal G_X[R_u(\mathcal G_X)-R_\gamma(\mathcal G_X)])\,du\biggr).
   \end{multline}
Here $\mathrm{Tr}(\mathcal G)$ is the trace of operator $\mathcal{G}$; $R_u(\mathcal G)$ is the resolvent operator defined by the formula
\begin{align*}
R_u(\mathcal G)=(I+2u\mathcal G)^{-1},
\end{align*}
and $\gamma=\gamma(\ep)$ satisfies the equation $\ep^2=\mathrm{Tr}(\mathcal{G}_X R_\gamma(\mathcal G_X))$.
\end{theorem}

In terms of eigenvalues $\mu_k$ of the covariance operator $\mathcal G_X$, this result reads as follows:

\begin{theorem}[{\cite[formula~(20)]{Syt1974}}]
    \label{thm:sytaya}
    If $\mu_k>0$ and $\sum\limits_{k=1}^\infty\mu_k<\infty$, then as $\ep\to0$
    \begin{align}
    \label{eq:sytaya}
        \mathbb{P}\Bigl\{\sum\limits_{k= 1}^\infty\mu_k\xi_k^2\leq \ep^2\Bigr\}
        \sim
        \bigl(-2\pi\gamma^2h''(\gamma)\bigr)^{-\frac 12}\exp\bigl(\gamma h'(\gamma)-h(\gamma)\bigr),
    \end{align}
    where $h(\gamma)=\frac12\sum\limits_{k=1}^\infty\log(1+2\mu_k\gamma)$ and $\gamma$ is uniquely determined by equation $\ep^2=h'(\gamma)$ for $\ep>0$ small enough.
\end{theorem}
As a particular case the formula~\eqref{eq:BB} for the Brownian bridge was also obtained. However, these theorems are not convenient to use for concrete computations and applications. For example, finding expressions for the function $h(\gamma)$ and determining the implicit relation $\gamma=\gamma(\ep)$ are the two difficulties that arise. The significant simplification was done only in 1997 by T.~Dunker, W.~Linde and M.~Lifshits~\cite{DLL1998}. We describe it in detail in Section~\ref{sec:secondwave}.

The results of~\cite{Syt1974} were not widely known. J.~Hoffmann-J{\o}rgensen~\cite{HoffJor1976} in 1976 obtained two-sided estimates for $L_2$-small ball probabilities. Later the result of Sytaya was rediscovered by I.A.~Ibragimov~\cite{Ibr1979}.\footnote{We notice that formula~(23) in~\cite{Ibr1979} contains several misprints. Apparently, this was first mentioned in~\cite[p.~745]{Fat2003}.} Cf. also \cite{NagSta1981}, \cite{Cox1982}, \cite{MayWolfZei1993}, \cite{DemZei1995}, \cite{Alb1996}. 

In 1984 V.M.~Zolotarev~\cite{Zol1986} suggested an explicit description of the small deviations in the case $\mu_k=\phi(k)$ with a decreasing and logarithmically convex function $\phi$ on $[1,\infty)$. His method was based on application of the Euler-Maclaurin formula to sums in the right hand side of~\eqref{eq:sytaya}. Unfortunately, no proofs were presented in \cite{Zol1986}, and, as was shown in~\cite{DLL1998}, this result is not valid without additional assumptions about the function $\phi$ (in particular, the final formula in~\cite[Example~2]{Zol1986} is not correct, see the end of Section~\ref{sec:secondwave} below).

In 1989 A.A. Borovkov and A.A. Mogul'skii~\cite{BorMog1989} investigated small deviations of several Gaussian processes in various norms. However, as noted in~\cite[Section~1]{Lif1999}, the result of Theorem~3 in~\cite{BorMog1989} contains an algebraic error.\footnote{On the other hand, the statement in~\cite[p.~745]{Fat2003} that this erroneous result was reprinted in~\cite[p.~206]{Lif1995-book} is in turn erroneous.}

One of the natural extensions of the problem under consideration is to find asymptotics~\eqref{eq:exact-small-ball-probability} for shifted small balls. It was shown in 
\cite[Theorem 1]{Syt1974} that if the shift $a$ belongs to the reproducing kernel Hilbert space\footnote{This space is often called the kernel of the distribution of the Gaussian vector $X$, see, e.g.,~\cite[Section~4.3]{Lif2012-book}.} (RKHS) of the Gaussian vector $X$, then
\begin{align}
\label{eq:shifted}
     \mathbb{P}\{\|X-a\|_{\mathfrak H}\leq\ep\}\sim \mathbb{P}\{\|X\|_{\mathfrak H}\leq\ep\}\cdot \exp\bigl(-\bigl\|\mathcal{G}_X^{-\frac 12}a\bigr\|_{\mathfrak{H}}^2/2\bigr),
\end{align}
where $\|\mathcal{G}_X^{-\frac 12}a\|_{\mathfrak{H}}$ is in fact the norm of $a$ in RKHS of $X$. See also some relations between the probabilities of centered and shifted balls  in~\cite{HoffJorSheDud1979} by~J.~Hoff\-mann-J{\o}rgensen, L.A.~Shepp and R.M.~Dudley, \cite{LindRos1994} by W.~Linde and J.~Rosinski, and~\cite{Der2003} by S.~Dereich. A generalized version of equivalence~\eqref{eq:shifted} for Gaussian Radon measures on locally convex vector spaces was proved by Ch.~Borell~\cite{Bor1977}.

Formula~\eqref{eq:shifted} was later generalized  by~J.~Kuelbs, W.V.~Li and W.~Linde~\cite{KueLiLin1994} and W.V.~Li, W.~Linde~\cite{LiLind1994}, who handled the asymptotics of probability 
\begin{align*}
    \mathbb{P}\{\|X-f(t)a\|_{\mathfrak{H}}\leq R(t)\},\qquad t\to\infty,
\end{align*}
for various combinations of $f$ and $R$, including the case $R(t)\to\infty$ as $t\to\infty$.

\section{Second wave}
\label{sec:secondwave}

The main difficulty in using Theorem~\ref{thm:sytaya} is that the explicit formulae for the eigenvalues of the covariance operators are rarely known. It was partially overcome by the celebrated \textit{\textbf{Wenbo Li comparison principle}}.\footnote{W.V. Li~\cite{Li1992} proved this statement provided $\sum|1-\mu_k/\tilde{\mu}_k|<\infty$ and conjectured that this assumption can be relaxed to the natural one $\prod\mu_k/\tilde{\mu}_k<\infty$. Later this conjecture was confirmed by F.~Gao, J.~Hannig, T.-Y.~Lee and F.~Torcaso~\cite{GaoHanLeeTor2004}.

W.~Linde~\cite{Lind1994-comparison} generalized W.V.~Li's assertion for the case of shifted balls.
In \cite{GaoHanTor2003} Theorem~\ref{thm:comparisonprinciple} was extended for the sums $\sum \mu_k |\xi_k|^p$, $p>0$, and even more general ones.}

\begin{theorem}[\cite{Li1992, GaoHanLeeTor2004}]
\label{thm:comparisonprinciple}
    Let $\xi_k$ be i.i.d. standard normal random variables, and let $(\mu_k)$ and $(\tilde{\mu}_k)$  be two positive non-increasing summable sequences such that  $\prod\mu_k/\tilde{\mu}_k<\infty$. Then
	\begin{align}
            \label{eq:Li0}
    	\mathbb{P}\Big\{\sum_{k=1}^{\infty} \tilde{\mu}_k \xi_k^2\leq\ep^2\Big\}
            \sim
    	\mathbb{P}\Big\{\sum_{k=1}^{\infty} \mu_k \xi_k^2\leq\ep^2\Big\}
    	\cdot \Big( \prod\limits_{k=1}^{\infty}\frac{\mu_k }{\tilde{\mu}_k} \Big)^{\frac 12},
    	\qquad\ep\to0.
	\end{align}	
\end{theorem}
\noindent So, if we know sufficiently sharp (and sufficiently simple!) asymptotic approximations $\tilde{\mu}_k$ for the eigenvalues $\mu_k$ of $\mathcal{G}_X$  then formula (\ref{eq:Li0}) provides the $L_2$-small ball asymptotics for $X$ up to a constant. However, to use this idea efficiently, we need explicit expressions of $L_2$-small ball asymptotics for ``model'' sequences $(\tilde{\mu}_k)$.
\medskip

The important step in the latter problem was made by M.A.~Lifshits in~\cite{Lif1997} who considered the small ball problem for more general series:
\begin{align}
    \label{eq:sum-general}
    \mathcal{S}:=\sum\limits_{k=1}^\infty \phi(k) Z_k,
\end{align}
where $(\phi(k))$ is a non-increasing summable sequence of positive numbers, and $Z_k$ are independent copies of a positive random variable $Z$ with finite variance and absolutely continuous distribution. The main restriction imposed in~\cite{Lif1997} on the distribution function $F(\cdot)$ of $Z$ is 
\begin{align}
\label{eq:F-ineq-Lifshits}
c_1 F(\tau)\leq F(b\tau)\leq c_2 F(\tau)
\end{align}
for some $b,c_1,c_2\in(0,1)$ and for all sufficiently small $\tau$. This assumption implies a polynomial (but not necessarily regular) lower tail behavior of the distribution. For the special case $Z_k=\xi_k^2$ the sum in~\eqref{eq:sum-general} corresponds to $L_2$-small ball probabilities~\eqref{eq:KL}.

Using the idea of R.A.~Davis and S.I.~Resnick~\cite{DavRes1991}, Lifshits expressed the exact small ball behavior of~\eqref{eq:sum-general} in terms of the Laplace transform of $\cal S$. If $Z$ has finite third moment, 
then a quantitative estimate of the remainder term was also given.
\medskip

The next important step was done by T.~Dunker, M.A.~Lifshits and W.~Linde in~\cite{DLL1998}. 
The result is based on the following assumption on the sequence $(\phi(k))$: 
\medskip

{\bf Condition DLL}. \textit{The sequence $(\phi(k))$ admits a  positive, logarithmically convex, twice differentiable and integrable extension on the interval $[1,+\infty)$.}\medskip

Under some additional assumptions\footnote{The results in~\cite{DLL1998} heavily depend on the ``condition {\bf I}\,'': total variation  $\mathbb{V}_{[0,\infty)} (\gamma f'(\gamma)/f(\gamma))$ is finite. This condition together with restriction~\eqref{eq:F-ineq-Lifshits} implies that the distribution function $F$ is regularly varying at zero with some index $\alpha<0$. Also this condition holds for the case~$Z_k=|\xi_k|^p$, $p>0$. Some relaxations are discussed in~\cite[Section 5]{DLL1998} and later papers of F.~Aurzada~\cite{Aur2007}, A.A.~Borovkov, P.S.~Ruzankin~\cite{BorRuz2008-1,BorRuz2008-2}, L.V.~Rozovsky~(see~\cite{Roz2017} and references therein).
} on the distribution function $F$ authors significantly simplified the expressions from~\cite{Lif1997} for the small ball behavior of $\cal S$.
By using Euler--Maclaurin's summation formula, they succedeed to express the small ball probabilities in terms of three integrals:
\begin{align*}
    I_0(\gamma):&=\int\limits_1^\infty (\log f)(\gamma\phi(t))\,dt,
    \\
    I_1(\gamma):&=\int\limits_1^\infty \gamma\phi(t) (\log f)'(\gamma\phi(t))\,dt,
    \\
    I_2(\gamma):&=\int\limits_1^\infty (\gamma\phi(t))^2 (\log f)''(\gamma\phi(t))\,dt.
\end{align*}
Here $f(\gamma)$ is the Laplace transform of $Z$,
\begin{align*}
f(\gamma)&:= \int\limits_0^\infty \exp(-\gamma\tau) \,dF(\tau).
\end{align*}

The main result (Theorem 3.1 in~\cite{DLL1998})  reads as follows:
\begin{align}
\label{eq:DLL-prob}
        \mathbb{P}\Bigl\{\mathcal{S}\leq r\Bigr\}\sim \sqrt{\frac{f(\gamma\phi(1))}{2\pi I_2(\gamma)}} \exp(I_0(\gamma)+\rho(\gamma)+\gamma r)),\qquad r\to0,
    \end{align}
    where $\gamma=\gamma(r)$ is any function satisfying
    \begin{align}
    \label{eq:DLL-ur}
        \lim\limits_{r\to0} \frac{I_1(\gamma)+\gamma r}{\sqrt{I_2(\gamma)}}=0,
    \end{align}
    and $\rho(\gamma)$ is a bounded function that represents the remainder terms in Euler--Maclaurin formula and can be written explicitly via infinite sums.

At a first glance, formula \eqref{eq:DLL-prob} and condition~\eqref{eq:DLL-ur} do not seem to be more explicit than Theorem~\ref{thm:sytaya}. Hovewer, for $Z_k=\xi_k^2$ (in this case $f(\gamma)=(1+2\gamma)^{-\frac 12}$) the result of \cite{DLL1998} turned out to be much more computationally tractable and became the base for derivation of the exact asymptopics for many ``model'' sequences of eigenvalues $(\phi(k))$.

In particular, the following asymptotics were obtained in~\cite[Section~4]{DLL1998}.

\begin{enumerate}
\item Let $p>1$. Then, as $\ep\to0$,
\begin{align}
\label{eq:DLL-polynomial}
\mathbb{P}\Big\{\sum_{k=1}^{\infty}k^{-p} \xi_k^2 \leq
\ep^2\Big\}
\sim 
{\cal C}\cdot \ep^{\gamma}
\exp\left(-\ D \ep^{-\frac
{2}{p-1}}\right),
\end{align}
where the constants $D$ and $\gamma$ depend on $p$ as follows: 
\begin{align*}
D=\frac {p-1}{2}\Big( \frac
{\pi}{p\sin\frac{\pi}{p}}\Big)^{\frac {p}{p-1}},\qquad \gamma = \frac {2-p}{2(p-1)},
\end{align*}
\noindent while the constant ${\cal C}$ is 
given by
the following expression:  
$$
{\cal C}=\frac
{(2\pi)^{\frac p4}\bigl(\sin \frac{\pi}{p}\bigr) ^{\frac {1 +
\gamma}{2}}} {(p-1)^{\frac 12} \bigl(\frac{\pi}{p}\bigr)^{1+\frac
{\gamma}{2}}}.
$$
Notice that this example was first considered in~\cite{Zol1986}.\footnote{However, as was first mentioned in~\cite[p.~743]{Fat2003}, the constants ${\cal C}$ and $D$ in~\cite[Example~1]{Zol1986} were calculated erroneously: they contain the Euler constant that does not appear in the correct answer. This error was reproduced in~\cite[formula~(3.2)]{Li1992}.}

\item  The second result shows that the general asymptotic formula in \cite{Zol1986} is not true. Namely, we have as $\ep\to0$
\begin{multline*}
\mathbb{P}\Bigl\{\sum\limits_{k=0}^\infty 
    \exp(-k)\,\xi_k^2\leq\ep^2\Bigr\}\\
\sim \frac{\exp\left(-\frac{\pi^2}{12}-\frac14(\log(\frac1{\ep^2}\log\frac1{\ep^2}))^2+\psi_0(\log(\frac1{\ep^2}\log\frac1{\ep^2}))\right)}{\pi^{\frac 12}\ep^{\frac 12}(\log\frac1{\ep^2})^{\frac 34}},
\end{multline*}
    where  $\psi_0$ is an explicitly given (though complicated) $1$-periodic and bounded function. It is shown in \cite{DLL1998} that $\psi_0$ is \textit{\textbf{non-constant}} while such term is absent in \cite[Example~2]{Zol1986}.\footnote{This error in~\cite{Zol1986} was also reproduced in~\cite[formula~(3.3)]{Li1992}.}

\end{enumerate}

\section{Exact asymptotics: Green Gaussian processes and beyond}
\label{sec:green}

\subsection{Green Gaussian processes and their properties}

As was mentioned above, the results of~\cite{Li1992}, \cite{GaoHanLeeTor2004} and \cite{DLL1998} allow to obtain the exact (at least up to a constant) small ball asymptotics for Gaussian random vectors from sufficiently sharp eigenvalues asymptotics of the corresponding covariance operators. The first significant progress in this direction was made for the special important class of Gaussian functions on the interval.

The \textit{\textbf{Green Gaussian process}} is a zero mean Gaussian process $X$ on the interval (say, $[0,1]$) such that its covariance function $G_X$ is the (generalized) Green function of an ordinary differential operator (ODO) on $[0,1]$ with proper boundary conditions. This class of processes is very important as it includes the Wiener process, the Brownian bridge, the Ornstein-Uhlenbeck process, their (multiply) integrated counterparts etc.
\medskip

First, we recall some definitions. 
Let ${\cal L}$ be an ODO given by the differential expression
\begin{align}
\label{eq:L}
{\cal L}u:=
(-1)^{\ell}\left(p_{\ell}u^{({\ell})}\right)^{({\ell})}+
\left(p_{\ell-1}u^{({\ell}-1)}\right)^{({\ell}-1)}+\dots+p_0u,
\end{align}
(here $p_j$, $j=0,\dots,\ell$ are functions on $[0,1]$, and
$p_{\ell}(t)>0$) and by $2\ell$ boundary conditions
\begin{align}
\label{eq:bc-non-separated}
U_{\nu}(u):= U_{\nu 0}(u)+U_{\nu 1}(u)=0,\qquad \nu=1,\dots,2\ell, 
\end{align}
where
\begin{align*}
U_{\nu 0}(u)&:=\alpha_{\nu}u^{(k_{\nu})}(0)+\sum\limits_{j=0}^{k_{\nu}-1}
\alpha_{\nu j}u^{(j)}(0),
\\
U_{\nu 1}(u)&:=\gamma_{\nu}u^{(k_{\nu})}(1)+\sum\limits_{j=0}^{k_{\nu}-1}
\gamma_{\nu j}u^{(j)}(1),
\end{align*}
and for any index $\nu$ at least one of coefficients $\alpha_{\nu}$ and 
$\gamma_{\nu}$ is not equal to zero.

For simplicity we assume $p_j\in W^j_{\infty}[0,1]$, $j=0,\dots,\ell$. 
Then the domain ${\cal D}({\cal L})$ consists of the functions 
$u\in W^{2\ell}_2[0,1]$ satisfying boundary conditions~\eqref{eq:bc-non-separated}.
\medskip

\begin{remark}
\label{rmk:kappa}
    It is well known, see, e.g., \cite[\S4]{Naim1969-book}, \cite[Chap.~XIX]{DunSch1971}, that the system of boundary conditions  can be reduced to the \textit{\textbf{normalized form}} by equivalent 
transformations. In what follows we always assume that this reduction is realized. This form is specified by the minimal sum of orders of all boundary conditions $\varkappa=\sum\limits_{\nu=1}^{2\ell} k_{\nu}$.
\end{remark}\medskip

The \textit{\textbf{Green function}} of the boundary value problem
\begin{align}
\label{eq:L-bvp}
{\cal L}u=\lambda u\quad {\rm on} 
\quad
[0,1], \qquad u\in {\cal D}({\cal L}),
\end{align}
is the function $G(t, s)$ such that it satisfies the equation ${\cal L}G(t,s) = \delta (s - t)$ in the sense of distributions and satisfies the boundary conditions~\eqref{eq:bc-non-separated}.\footnote{Notice that if $G$ is the covariance function of a random process then $G(t, s)\equiv G(s,t)$, and therefore the problem~\eqref{eq:L-bvp} is always self-adjoint.} The existence of the Green function is equivalent to the invertibility of the operator ${\cal L}$ with given boundary conditions, and $G(t, s)$ is the kernel of the integral operator~${\cal L}^{-1}$. 

If the problem~\eqref{eq:L-bvp} has a zero eigenvalue corresponding to the eigenfunction $\varphi_0$ (without loss of generality, it can be assumed to be normalized in $L_2[0, 1]$), then the Green function obviously does not exist. If $\varphi_0$ is unique up to a constant multiplier,  then the function $G(t, s)$ is called the \textit{\textbf{generalized Green function}} if it satisfies the equation ${\cal L}G(t, s) = \delta(t - s) - \varphi_0(t)\varphi_0(s)$ in the sense of distributions, subject to the boundary conditions and the orthogonality condition 
\begin{align*}
\int\limits_0^1 G(t,s)\,\varphi_0(s)\,ds=0,\qquad 0 \leq t \leq 1.
\end{align*}
The generalized Green function is the kernel of the integral operator which is inverse to ${\cal L}$ on the subspace of functions orthogonal to $\varphi_0$ in $L_2[0, 1]$. In a similar way, one can consider the case of multiple zero eigenvalue.

\medskip 

Thus, (non-zero) eigenvalues $\mu_k$ of the covariance operator of a Green Gaussian process are inverse to the (non-zero) eigenvalues $\lambda_k$ of the corresponding boundary value problem~\eqref{eq:L-bvp}: $\mu_k=\lambda_k^{-1}$. So, to obtain rather good asymptotics of eigenvalues $\mu_k$ one can use the powerful methods of spectral theory of ODOs, originated from the classical works of G.~Birkhoff~\cite{Birk1908-1}, \cite{Birk1908-2} and
J.D.~Tamarkin~\cite{Tam1912}, \cite{Tam1928}. 

\medskip

Notice that classical operations on the random processes -- integration and centering -- transform a Green Gaussian process to a Green one. It is easy that if $G(t,s)$ is the covariance function of a random process $X$ then the covariance functions of the integrated and the centered (\textit{\textbf{demeaned}}) process 
\begin{align}
\label{eq:integrated-centered-process}
X_1(t)=\int\limits_0^t X(s)\,ds, \qquad \overline{X}(t)=X(t)-\int\limits_0^1 X(s)\,ds
\end{align}
are, respectively,
\begin{align}
\label{eq:G-integrated}
G_1(t,s)&=\int\limits_0^t\int\limits_0^sG(x,y)\, 
dydx,
\\
\label{eq:G-centered}
{\overline G}(t,s)&=G(t,s)-\int\limits_0^1G(t,y)\,dy-\int\limits_0^1G(x,s)\,dx+\int\limits_0^1\int\limits_0^1G(x,y)\,dydx.
\end{align}

\begin{theorem}
\label{thm:Green-to-Green}
{\bf 1} {\rm ({\cite[Theorem~2.1]{NazNik2004}}).} Let the kernel $G(t,s)$ be the Green
function for the boundary value problem~\eqref{eq:L-bvp}. 
Then the integrated kernel~\eqref{eq:G-integrated} is the Green function for the  boundary value problem 
\begin{align}
\label{eq:L1-integrated-bvp}
{\cal L}_1u:=-({\cal L}u')'=\lambda u\quad {\rm
on}\quad [0,1], \qquad u\in {\cal D}({\cal L}_1),
\end{align}
where the domain ${\cal D}({\cal L}_1)$ consists of
functions $u\in W_2^{2\ell+2}[0,1]$ satisfying the
boundary conditions
\begin{align*}
u(0)=0; \qquad u'\in
{\cal D}({\cal L}); \qquad
({\cal L}u')(1)=0.
\end{align*}

{\bf 2} {\rm  ({\cite[Theorem~3.1]{Naz2009-non-separated}}).} Let the boundary value problem~\eqref{eq:L-bvp} have a zero eigenvalue 
with constant eigenfunction $\varphi_0(t)\equiv1$, and let the kernel $G(t,s)$ 
be the generalized Green function of the problem~\eqref{eq:L-bvp}. Then the integrated 
kernel~\eqref{eq:G-integrated} is the (conventional) Green function of the boundary value problem~\eqref{eq:L1-integrated-bvp} where the domain ${\cal D}({\cal L}_1)$ consists of functions
$u\in W_2^{2\ell+2}[0,1]$ satisfying the boundary conditions
\begin{align}
\label{eq:bc-integrated}
u(0)=0;\qquad u(1)=0;\qquad u'\in {\cal D}({\cal L}).
\end{align}

{\bf 3} {\rm (\cite[Theorem 1]{NazPet2023}\footnote{For the problem~\eqref{eq:L1-integrated-bvp}--\eqref{eq:bc-integrated} this statement was proved in~\cite{Naz2009-non-separated}.}).} Let the kernel $G(t,s)$ be the Green function of the problem~\eqref{eq:L-bvp}, and let the corresponding differential expression ${\cal L}$ have no zero order term\footnote{The case $p_0\not\equiv0$ is much more complicated. In this case the centered process~$\overline{X}$ is in general not a Green Gaussian process.} ($p_0\equiv0$).   
Then the centered kernel~\eqref{eq:G-centered} is the generalized Green function for the boundary value problem $\overline{\cal L}u:={\cal L}u=\lambda u$ with explicitly given boundary conditions.
\end{theorem}

We provide three examples of transformations of Green Gaussian processes generating somewhat more complicated boundary value problems: 
\begin{enumerate}
    \item Multiplication by a deterministic function $\psi$;
    \item The so-called \textit{\textbf{online-centering}}~\cite{KlepLeBre2002}:
\begin{align*}
\widehat{X}(t)=X(t)-\frac1t\int\limits_{0}^t X(s)\,ds;
\end{align*}
    \item Linear combination of two Green Gaussian processes.
\end{enumerate}

 \begin{theorem}
{\bf 1}~\mbox{\rm(\cite[Lemma~2.1]{NazPus2009})}.\footnote{For $X=W$ and $X=B$ this fact was obtained earlier in \cite[Theorems~1.1 and 1.2]{DehMar2003}. Some interesting examples can be found in~\cite{Pyck2021}.}
Let $X(t)$ be a Green Gaussian process, corresponding to the boundary value problem~\eqref{eq:L-bvp}, and let $\psi\in W^{\ell}_\infty[0,1]$, $\psi>0$ on $(0,1)$. Then the covariance function of the process $\psi(t)X(t)$ is the Green function of the problem
\begin{align*}
{\cal L}_{\psi}u:=\psi^{-1}{\cal L}(\psi^{-1} u)=\lambda u\quad \mbox 
{\rm on}\quad [0,1], \qquad \psi^{-1}u\in {\cal D}({\cal L}).
\end{align*}

{\bf 2}~\mbox{\rm(\cite[Theorem~2]{NazNik2022-spect-equiv})}. Let $X(t)$ be a Green Gaussian process, corresponding to the boundary value problem~\eqref{eq:L-bvp}. Then the covariance function of the online-centered integrated process $\widehat{X_1}(t)$ is the Green function of the problem
\begin{align*}
\widehat{{\cal L}_1}u:=t(t^{-1}{\cal L}(t^{-1}(tu)'))'=\lambda u\quad \mbox 
{\rm on}\quad [0,1], \qquad u\in {\cal D}(\widehat{{\cal L}_1}),
\end{align*}
\noindent where the domain ${\cal D}(\widehat{{\cal L}_1})$ is defined by $2\ell+2$ boundary conditions
\begin{align*}
u(0)=0,\qquad t^{-1}(tu)'\in{\cal D}({\cal L}),\qquad {\cal L}(t^{-1}(tu)')(1)=0.
\end{align*}

{\bf 3}~\mbox{\rm(\cite[Section~2]{NazNik2018-mixed})}. Let $X_1(t)$ and $X_2(t)$ be independent Green Gaussian processes, corresponding to the boundary value problems of the form~\eqref{eq:L-bvp} with the same operator ${\cal L}$ and different boundary conditions. Then for any $\alpha_1,\alpha_2>0$ such that $\alpha_1^2+\alpha_2^2=1$, the \textit{\textbf{mixed}} process
$$
X(t)=\alpha_1 X_1(t)+\alpha_2 X_2(t)
$$
is a Green Gaussian process, corresponding to the boundary value problem of the form~\eqref{eq:L-bvp} with the same operator ${\cal L}$ subject to some (in general, more complicated) boundary conditions.
\end{theorem}

\subsection{Exact $L_2$-small ball asymptotics\\ for the Green Gaussian processes}

Apparently, the papers~\cite{BegNikOrs2005} by 
L.~Beghin, Ya.Yu.~Nikitin and E.~Orsingher and~\cite{DehMar2003} by P.~Deheuvels and G.~Martynov were the first where the eigenvalues asymptotics for covariance operators was obtained by the passage to the corresponding boundary value problem. The small deviation asymptotics up to a constant was derived for several concrete Green Gaussian processes, including  some integrated 
ones~\cite{BegNikOrs2005} and weighted ones~\cite{DehMar2003}.

To formulate the general result, first, we describe the asymptotic behavior of the eigenvalues of the boundary value problem~\eqref{eq:L-bvp}.

It is well known, see, e.g., \cite[\S4]{Naim1969-book}, \cite[Chap.~XIX]{DunSch1971}, that eigenvalues $\lambda_k$ of Birkhoff-regular (in particular, self-adjoint) boundary value problem for ODOs with smooth coefficients can be expanded into asymptotic series in powers of  $k$ (analogous results under more general hypotheses, as well as some additional references, can be found in~\cite{Shka1982}, \cite{Shka1983}). Taking into account the Li comparison principle (Theorem~\ref{thm:comparisonprinciple}), it is easy to show that the $L_2$-small ball asymptotics up to a constant for the Green Gaussian process requires just two-term spectral asymptotics for corresponding boundary value problem (with the remainder estimate). This asymptotics is given by the following assertion.

\begin{theorem}[{\cite[\S4, Theorem 2]{Naim1969-book}}] Let $p_{\ell}\equiv1$. Then the eigenvalues of the problem~\eqref{eq:L-bvp} counted according to their 
multiplicities can be split into two subsequences $(\lambda'_j)$, $(\lambda''_j)$, 
$j\in\mathbb{N}$, such that, as $j\to\infty$, 
\begin{align}
\label{eq:lambda-two-term}
\lambda'_j=\left(2\pi j+\rho'+O(j^{-\frac12})\right)^{2\ell}, 
\qquad
\lambda''_j=\left(2\pi j+\rho''+O(j^{-\frac12})\right)^{2\ell}.
\end{align}
\end{theorem}

The first term of this expansion is completely determined by the principal coefficient of the operator,\footnote{In general case the problem \eqref{eq:L-bvp} can be reduced to the case  
$p_{\ell}\equiv1$ by the independent variable transform, see \cite[\S4]{Naim1969-book}. 
The expressions in brackets in~\eqref{eq:lambda-two-term} should be divided by
$\int\limits_0^1 p_{\ell}^{-\frac 1{2\ell}}(t)\,dt$.} whereas the formulae for other terms (in particular, for $\rho'$ and $\rho''$) are rather complicated. However, it turns out that the exact $L_2$-small ball asymptotics (up to a constant) for the Green Gaussian processes depends only on the sum $\rho'+\rho''$, which is completely determined by~$\varkappa$, the sum of orders of boundary conditions, see Remark~\ref{rmk:kappa}. So we arrive at the following result.

\begin{theorem}[{\cite[Theorem~7.2]{NazNik2004}}; {\cite[Theorem~1.2]{Naz2009-non-separated}}]
\label{thm:small-ball-Green}
 Let the covariance function $G_X(t,s)$ of a zero mean 
Gaussian process $X(t)$, $0\le t \le 1$, be the Green function of the boundary value problem~\eqref{eq:L-bvp} generated by a 
differential expression~\eqref{eq:L} and by boundary conditions~\eqref{eq:bc-non-separated}. Let 
$\varkappa<2\ell^2$. Then, as $\ep\to0$,
\begin{align}
\label{eq:P-non-separated}
\mathbb{P}\{{}\|X\|_{2,[0,1]}\le \ep\}\sim {\cal C}(X)\cdot 
\ep^{\gamma} \exp\Big(-D_\ell\ep^{-\frac {2}{2\ell-1}} \Big).
\end{align}
Here we denote 
\begin{align}
\label{eq:gamma}
D_\ell=\ \frac {2\ell-1}{2}
\left(\frac {\theta_{\ell}}{2\ell\sin\frac {\pi}{2\ell}}\right)
^{\frac {2\ell}{2\ell-1}},
\quad
\gamma=-\ell+\frac {\varkappa+1}{2\ell-1},
\quad 
\theta_{\ell}=\int\limits_0^1p_{\ell}^{-\frac 1{2\ell}}(t)\,dt,
\end{align}
and the constant ${\cal C}(X)$ is given by 
\begin{align}
\label{eq:C-constant}
{\cal C}(X)=
C_{\rm dist}(X)\cdot\frac {(2\pi)^{\frac\ell2}
\bigl(\frac\pi{\theta_{\ell}}\bigr)
^{\ell\gamma}\left(\sin \frac {\pi}{2\ell}\right)
^{\frac {1 +\gamma}{2}}} 
{(2\ell-1)^{\frac12} 
\left(\frac {\pi}{2\ell}\right)^{1+\frac {\gamma}2}
{\Gamma\vphantom)}^{\ell}\left(\ell-\frac {\varkappa}{2\ell}\right)},
\end{align}
where $C_{\rm dist}(X)$ is the so-called {\bf\textit{ distortion constant}} (cf. formula~\eqref{eq:Li0})
\begin{align}
\label{eq:C-dist}
C_{\rm dist}(X)
:=
\prod_{k=1}^{\infty}
\frac {\lambda_k^{\frac12}}{\bigl(\frac\pi{\theta_{\ell}} \cdot 
\left[k+\ell-1-\frac {\varkappa}{2\ell}\right]\bigr)^{\ell}},
\end{align}

\noindent and $\lambda_k$ are the eigenvalues of the problem~\eqref{eq:L-bvp} taken in non-decreasing order counted with their multiplicities.
\end{theorem}

\begin{remark}
\label{rmk:zero-eigen-1}
If $G_X(t,s)$ is the generalized Green function of the problem~\eqref{eq:L-bvp}
then formula~\eqref{eq:P-non-separated} 
should be corrected. Namely, if, say, one of eigenvalues $\lambda_k$ vanishes,  the exponent $\gamma$ in~\eqref{eq:gamma} increases by $\frac {2\ell}{2\ell-1}$, and 
$\lambda_k$ in the product \eqref{eq:C-dist} should be properly renumbered.

On the other hand, if, say, $\varkappa=2\ell^2$ then formulae \eqref{eq:C-constant} and \eqref{eq:C-dist} cannot be applied directly. The correct result can be obtained if we substitute $\varkappa\to\varkappa-\epsilon$ and then push $\epsilon\to0$ in \eqref{eq:C-constant}.
\end{remark}

The first step towards a general result was made in~\cite{GaoHanTor2003-integrated} by F.~Gao, J.~Hannig and F.~Torcaso. They considered $m$-times integrated Wiener processes $W_m^{[\beta_1,\ldots,\beta_m]}(t)$. Here and later on, for a random process $X(t)$ on $[0,1]$ we denote   
\begin{align}
\label{eq:m-times-interated-X}
X_m^{[\beta_1,\ldots,\beta_m]}(t)
=
(-1)^{\beta_1+\ldots+\beta_m}\int\limits_{\beta_m}^t\int\limits_{\beta_{m-1}}^{t_{m-1}}\ldots\int\limits_{\beta_1}^{t_1}X(s)\,dsdt_{1}\ldots\,dt_{m-1}
\end{align}
(here $0 \leq t \leq 1$, and any $\beta_j$ equals either zero or one;  the usual $m$-times integrated process corresponds to $\beta_1=\dots=\beta_m=0$).

From Theorem~\ref{thm:Green-to-Green} it follows that the covariance function of $W_m^{[\beta_1,\ldots,\beta_m]}$ is the Green function of the boundary value problem
\begin{align*}
(-1)^{m+1}u^{(2m+2)}=\lambda u\quad {\rm
on}\quad [0,1], 
\end{align*}
with boundary conditions depending on parameters $\beta_j$, $j=1,\dots,m$. The special case $W_m^{[1,0,1,\ldots]}$ is called in~\cite{GaoHanTor2003-integrated} the Euler-integrated Brownian motion, since its covariance function can be expressed in terms of Euler polynomials. 

By particular fine analysis, authors of~\cite{GaoHanTor2003-integrated}
obtained the formula
\begin{align}
\label{eq:small-ball-integrated-w-gao}
\mathbb{P}\{{}\|W_m^{[\beta_1,\ldots,\beta_m]}\|_{2,[0,1]}\le \ep\}
\sim 
C\cdot 
\ep^{\frac{1-k_0(2m+2)}{2m+1}} 
\exp\Big(- D_m\ep^{-\frac {2}{2m+1}} \Big).
\end{align}
Here $D_m$ is equal to the corresponding coefficient in~\eqref{eq:P-non-separated} for $\ell=m+1$ and $k_0$ is an \textit{unknown} integer. It was also conjectured in~\cite{GaoHanTor2003-integrated} that $k_0=0$ for all $m\in\mathbb{N}$ and for any choice of $\beta_j$ (for $m\leq 10$ this was checked numerically). This conjecture was verified in~\cite{GaoHanLeeTor2004}. 

Simultaneously, A.I. Nazarov and Ya.Yu. Nikitin~\cite{NazNik2004} obtained the result of Theorem~\ref{thm:small-ball-Green} for arbitrary Green Gaussian process under assumption that the boundary conditions~\eqref{eq:bc-non-separated} 
are \textit{\textbf{separated}}:
\begin{align}
\label{eq:bc-separated}
U_{\nu 0}(u)=0,\quad \nu=1,\dots,\ell; \qquad
U_{\nu 1}(u)=0,\quad \nu=\ell+1,\dots,2\ell. 
\end{align}
For general self-adjoint boundary value problems, Theorem~\ref{thm:small-ball-Green} was established later in~\cite{Naz2009-non-separated}.

Theorem~\ref{thm:small-ball-Green} provides the exact $L_2$-small ball asymptotics up to a constant. Methods of derivation of this constant were proposed simultaneously and independently in~\cite{Naz2003-sharp} and~\cite{GaoHanLeeTor2003-Hadamard}. Although slightly different, both methods are based on complex variables techniques (the Hadamard factorization and the Jensen theorem).

\begin{remark}
As an auxiliary result, the following generalization of~\eqref{eq:DLL-polynomial} was obtained in~\cite[Theorem~6.2]{NazNik2004}:\footnote{For $p=2$ this formula was established earlier in \cite{Li1992}.} 

Let $\th>0$, $\delta>-1$, $p>1$. Then, as $\ep\to0$,
\begin{align}
\label{eq:P-model-naznik}
\mathbb{P}\Big\{\sum_{k=1}^{\infty}(\th (k+\delta))^{-p} \xi_k^2 \leq
\ep^2\Big\}
\sim 
{\cal C}(\th,p,\delta)\cdot \ep^{\gamma}
\exp\Big(-D\ep^{-\frac {2}{p-1}} \Big),
\end{align}
where 
\begin{align}
\label{eq:D-gamma}
D=\frac {p-1}{2}\Big( \frac
{\pi}{\th p\sin\frac{\pi}{p}}\Big)^{\frac {p}{p-1}},\qquad \gamma=\frac{2-p-2\delta p}{2(p-1)},
\end{align}
while
$$
{\cal C}(\th,d,\delta) =\frac
{(2\pi)^{\frac{p}4}\th^{\frac{d\gamma}2}\bigl(\sin \frac{\pi}p\bigr) ^{\frac {1 +
\gamma}{2}}} {(p-1)^{\frac12} \bigl(\frac{\pi}p\bigr)^{1+\frac
{\gamma}{2}}\Gamma^{\frac{p}2}(1+\delta)}.
$$

Formula~\eqref{eq:P-model-naznik}, together with Theorem \ref{thm:comparisonprinciple} and the results of \cite{Naz2003-sharp}, gave the opportunity to derive (\cite{NikKha2006}, \cite{Roz2021}) the exact $L_2$-small ball asymptotics for Gaussian processes satisfying the following condition on the eigenvalues $\mu_k$ of the covariance operator:
\begin{align*}
\mu_k=\frac {\big(P_1(k)\big)^{\nu_1}}{\big(P_2(k)\big)^{\nu_2}},
\end{align*}
where $P_1$ and $P_2$ are polynomials, $\nu_1,\nu_2>0$.\footnote{A particular case $\mu_k=\frac 1{(k-a)(k-b)}$ was considered earlier in~\cite[Theorem~3.7]{Fat2003}.} Several examples of such processes can be found in~\cite{Pyck2003-phd} (see also \cite{Pyck2001}, \cite{Pyck2003-multivariate}).
\end{remark}

Now we give an important corollary from Theorem~\ref{thm:small-ball-Green} and Part~1 of Theorem~\ref{thm:Green-to-Green}.

\begin{corollary}[{\cite[Theorem~4.1]{NazNik2004}}]
\label{cor:integrated-constant}
Given Green Gaussian process $X$ on $[0,1]$, the parameter $\varkappa$ (the sum of orders of boundary conditions) for a boundary value problem~\eqref{eq:L-bvp} corresponding to any $m$-times integrated process $X_m^{[\beta_1,\ldots,\beta_m]}$
does not depend on $\beta_j$, $j=1,\dots,m$. Thus, the $L_2$-small ball asymptotics for various processes $X_m^{[\beta_1,\ldots,\beta_m]}$ can differ only by a constant.
\end{corollary}

This corollary generates a natural question: which choices of parameters $\beta_j$ give 
extremal (maximal/minimal) constant in $L_2$-small ball asymptotics among all $m$-times integrated processes?

For $X=W$ the answer was given independently in \cite[Remark~2]{GaoHanLeeTor2004} and \cite[Theorem~4.1]{Naz2003-sharp}. It turns out that the usual integrated Wiener process  $W_m^{[0,0,0,\ldots]}$ has the biggest multiplicative constant, while the Euler-integrated process $W_m^{[1,0,1,\ldots]}$ has the smallest one. The same extremal properties of the processes $X_m^{[0,0,0,\ldots]}$ and $X_m^{[1,0,1,\ldots]}$ for the Green Gaussian processes $X$ under some symmetry assumptions\footnote{In particular, these assumptions are fulfilled for a symmetric process, for instance, for the Brownian bridge and the Ornstein--Uhlenbeck process.} were obtained in~\cite[Proposition~4.4]{Naz2003-sharp},~\cite{NazNaz2004}.\medskip

Methods established in~\cite{NazNik2004},~\cite{GaoHanLeeTor2004},~\cite{Naz2009-non-separated} and~\cite{Naz2003-sharp},~\cite{GaoHanLeeTor2003-Hadamard} were widely used to obtain exact small ball asymptotics to many concrete Green Gaussian processes.  Several examples were considered even in the pioneer papers~\cite{Naz2003-sharp}, \cite{GaoHanLeeTor2003-Hadamard}. Much more results can be found in~\cite{NikOrs2004},~\cite{Deh2006},~\cite{Naz2009-non-separated},~\cite{Pus2010-bogoliubov},~\cite{Pus2011-Matern},~\cite{BarIgl2011}.\footnote{The thesis~\cite{Jin2014-phd}
 is of particular interest. Besides several known examples, the Slepian process $S(t)=W(t+a)-W(t)$ on $[0,1]$ was considered there. It turns out that for $a<1$ eigenproblem for corresponding covariance operator is equivalent to the boundary value problem for a \textit{\textbf{system}} of ordinary differential equations. The analysis of this boundary value problem is much more difficult than in the case of a single equation; in~\cite[Sections 2.2.10 and 3.3.7]{Jin2014-phd} it is managed only for $\frac 12\le a<1$.} 

\medskip

The $L_2$-small ball asymptotics of the weighted Green Gaussian processes are of particular interest. The first results on exact asymptotics were obtained in
\cite{DehMar2003}, \cite{Naz2003-sharp}, \cite{GaoHanLeeTor2003-Hadamard}. 

Numerous examples were considered in~\cite{NazPus2009}. Also it was noticed in this paper that if the weight is sufficiently smooth and \textit{\textbf{non-degenerate}} (i.e. bounded and bounded away from zero)  
then the formula~\eqref{eq:P-non-separated} holds with exponent $\gamma$ independent of the weight.\footnote{In contrast, the coefficient $D_\ell$ depends on the weight.} However, if the weight degenerates
(vanishes or blows up) at least at one point of the segment, generally speaking, 
this is not the case. 

In fact, the following comparison theorem holds for the Green Gaussian processes with non-degenerate weights. For simplicity we give it here for the case of separated boundary conditions, a more general result can be found in~\cite[Corollary 1]{NazPus2014-comparison}. 

\begin{theorem}[{\cite[Theorem~2]{NazPus2014-comparison}}]\!\!\footnote{For some concrete processes this result was established earlier by Ya.Yu.~Nikitin and R.S.~Pusev~\cite{NikPus2013}.}
\label{spectral}
Let $X(t)$ be a Green Gaussian process on $[0,1]$. Suppose that the boundary conditions of corresponding boundary value problem are separated (see~\eqref{eq:bc-separated}).

Suppose that the weight functions $\psi_1, \psi_2\in W_\infty^{\ell}[0,1]$ are bounded away from zero, and
\begin{equation*}
\int\limits_0^1\psi_1^{\frac1{\ell}}(t)\,dt=\int\limits_0^1\psi_2^{\frac1{\ell}}(t)\,dt.
\end{equation*}
Denote by $\varkappa_0$ and $\varkappa_1$ sums of orders of boundary conditions at zero and one, respectively:
$\varkappa_0=\sum\limits_{\nu=1}^{\ell} k_{\nu}$, $\varkappa_1=\sum\limits_{\nu=\ell+1}^{2\ell} k_{\nu}$.
Then, as $\ep\to0$,
\begin{multline*}
\mathbb{P}\bigl\{\|\psi_1X\|_{2,[0,1]}\leq\ep\bigr\}\\
\sim\mathbb{P}\bigl\{\|\psi_2 X\|_{2,[0,1]}\leq\ep\bigr\}\cdot 
\Big(\frac{\psi_2(0)}{\psi_1(0)}\Big)
^{-\frac{\ell}2+\frac14+\frac{\varkappa_0}{2\ell}}
\Big(\frac{\psi_2(1)}{\psi_1(1)}\Big)
^{-\frac{\ell}2+\frac14+\frac{\varkappa_1}{2\ell}}.
\end{multline*}
\end{theorem}

\subsection{Finite-dimensional perturbations\\ of Gaussian random vectors}

A natural class of Gaussian random vectors that often can be treated analytically is the class of finite-dimensional perturbations of the Gaussian random vector $X$ for which we already know the exact small deviation asymptotics.\footnote{Notice that the \textit{\textbf{logarithmic}} asymptotics does not change under finite-dimensional perturbations, see Section~\ref{sec:logarithm}.} In this Subsection we consider two types of such perturbations arising in applications.

\subsubsection{Perturbations from the kernel of original Gaussian random vector}

The first example of such perturbations was considered by P.~Deheuvels~\cite{Deh2007}. He introduced the following process:
\begin{align*}
Y_K(t)=B(t)-6K\, t(1-t)\int\limits_0^1 B(s)\,ds,\qquad 0\leq t\leq 1,
\end{align*}
where $B(t)$ is the standard Brownian bridge while $\alpha\in\mathbb{R}$ is a constant. He showed that the distributional equality $Y_K(t)\stackrel{d}{=}Y_{2-K}(t)$ holds, and obtained the explicit Karhunen--Lo{\`e}ve expansion for the process $Y_1(t)$. 

A general construction of one-dimensional perturbations in this class was introduced by A.I.~Nazarov~\cite{Naz2010-onedimpert}.

Suppose that $\mathcal{O}\subset\mathbb{R}^d$ is a bounded domain, and let $X(x)$ be a Gaussian random function on $\overline{\mathcal{O}}$. Consider the following family of ``perturbed'' Gaussian random functions: 
\begin{align}
\label{eq:onedimpert}
X_{\varphi,\alpha}(x)=X(x)-\alpha\psi(x)\int\limits_{\mathcal{O}} X(u)\varphi(u)\, du.
\end{align}
Here $\alpha\in\mathbb{R}$, $\varphi$ is a locally summable function in $\mathcal{O}$, $\psi=\mathcal{G}_X\varphi$, and
\begin{align}
\label{eq:q}
q:=\int\limits_{\mathcal O}\psi(x)\varphi(x)\,dx=\iint\limits_{\mathcal{O}\times\mathcal{O}} G_X(x,y) \varphi(x)\varphi(y)\,dxdy<\infty.
\end{align}
The latter condition just means that $\psi$ belongs to the kernel of the distribution of $X$.\footnote{In fact, the results in \cite{Naz2010-onedimpert} hold in a more general setting, where $X$ is a Gaussian vector in the Hilbert space, $\varphi$ is a measurable linear functional of $X$, and $\psi=\mathcal{G}_X\varphi$.
}

\begin{remark}
It is easy to see that the Deheuvels process $Y_K$ is a particular case of $X_{\varphi,\alpha}$ for $X=B$, $\varphi\equiv1$ and $\alpha=12K$.
\end{remark}

Direct calculation shows that the covariance function $G_{\varphi,\alpha}$ of $X_{\varphi,\alpha}$ is a one-dimensional perturbation of the covariance function $G$, namely,
\begin{align*}
G_{\varphi,\alpha}(x,y)=G_X(x,y)+Q\psi(x)\psi(y),
\end{align*}
where $Q=
q\bigl(\alpha-\frac 1q\bigr)^2-\frac 1q$. The following statements can also be easily checked.

\begin{lemma}[{\cite[Corollaries~1,~2]{Naz2010-onedimpert}}]
\begin{itemize}
    \item[{\bf 1.}] For the random functions~\eqref{eq:onedimpert}, the equality
	$X_{\varphi,\alpha}(x)\stackrel{d}{=}X_{\varphi,\frac2q-\alpha}(x)$
holds.
\item[{\bf 2.}] Let $\widehat\alpha=\frac 1q$. Then
\begin{itemize}
    \item[a)]
The identity
$\int\limits_{\cal O}{\cal X}_{\varphi,\widehat\alpha}(x)\varphi(x)\,dx=0$ holds true a.s.
\item[b)]
The process ${\cal X}_{\varphi,\widehat\alpha}(x)$ and the random variable
$\int\limits_{\cal O} X(u)\varphi(u)du$ are independent.
\end{itemize}
\end{itemize}
\end{lemma}

If $\alpha=\frac 1q$ then the random function \eqref{eq:onedimpert} is called the \textit{\textbf{critical perturbation}} of $X$. Otherwise it is called the non-critical perturbation.

It turns out that in the non-critical case the $L_2$-small ball asymptotics for ${\cal X}_{\varphi,\alpha}(x)$ coincides with that of $X$ up to a constant.

\begin{theorem}[{\cite[Theorem~1]{Naz2010-onedimpert}}]
\label{thm:onedimpert-non-critical}	
Let a function $\varphi$ satisfy~\eqref{eq:q}. 
If $\alpha\not=\frac1q$ then, as $\ep\to0$,
	\begin{align*}
		\mathbb{P}\left\{\|X_{\varphi,\alpha}\|_{2,\mathcal{O}}\leq\ep\right\}\sim |1-\alpha q|^{-1}\cdot\mathbb{P}\left\{\|X\|_{2,\mathcal{O}}\leq\ep\right\}.
	\end{align*}
\end{theorem}
For the critical perturbations, a general result is essentially more complicated, see~{\cite[Theorem~2]{Naz2010-onedimpert}}. We restrict ourselves to the case where $X$ is a Green Gaussian process on $[0,1]$. We underline that the requirement on $\varphi$ here is stronger than in Theorem~\ref{thm:onedimpert-non-critical}.

\begin{theorem}[{\cite[Theorem~3]{Naz2010-onedimpert}}]
\label{thm:onedimpert-Green}
Let the covariance function $G_X$ be the Green function for the boundary value problem~\eqref{eq:L-bvp}, where ${\cal L}$ is an ODO~\eqref{eq:L} with $p_\ell\equiv 1$. Suppose that $\hat{\alpha}=\frac 1q$ and $\varphi\in L_2[0,1]$. Then, as $\ep\to0$,
\begin{align*}
\mathbb{P}\{\|X_{\varphi,\hat{\alpha}}\|_{2,[0,1]}\leq\ep\}
\sim 
\frac{\sqrt{q}}{\|\varphi\|_{2,[0,1]}}\cdot \Bigl(2\ell\sin\bigl(\frac\pi{2\ell}\bigr)\ep^2\Bigr)^{-\frac{\ell}{2\ell-1}}\cdot \mathbb{P}\{\|X\|_{2,[0,1]}\leq\ep\}.
\end{align*}
\end{theorem}

Also in~\cite{Naz2010-onedimpert} an algorithm of derivation of Karhunen--Lo{\`e}ve expansion for the process
${\cal X}_{\varphi,\alpha}$ was given\footnote{In a particular case this algorithm was invented by M.~Kac, J.~Kiefer and J.~Wolfowitz in~\cite{KKW1955}. Some examples can be also found in a recent paper~\cite{BarLov2018}.} provided we know explicitly 
the fundamental system of solutions to the equation ${\cal L}u-\lambda u=0$ for $\lambda\in\mathbb{R}$. For instance, this is the case if $X$ is the Brownian bridge, which is very important in applications (see below).

The results of~\cite{Naz2010-onedimpert} were generalized for the case of multi-dimensional perturbations of Gaussian random functions by Yu.P.~Petrova~\cite{Pet2021}.\footnote{The problem of $L_2$-small ball asymptotics for some concrete examples of such perturbations of the Brownian bridge appearing in Statistics in a regression context was studied earlier in~\cite{KirNik2014}. Notice, however, that the power term in~\cite[Theorem~3.1]{KirNik2014} was calculated erroneously.}
\medskip

Finite-dimensional perturbations of the Wiener process or the Brownian bridge often appear in Statistics. Let us consider an $\omega^2$-type goodness-of-fit test with parameters estimated from the sample. Corresponding limiting process in this case (the so-called Durbin-type process, see~\cite{Dur1973}) is just a finite-dimensional perturbation\footnote{Namely, critical perturbation, as it was shown in~\cite[Theorem~4]{Pet2021}.} of the Brownian bridge. 

Unfortunately, in many important cases the perturbation functions $\varphi$ for the Durbin-type processes do not belong to $L_2[0,1]$, and thus Theorem~\ref{thm:onedimpert-Green} is not applicable. For instance, this is the case for the Kac--Kiefer--Wolfowitz processes~\cite{KKW1955} appearing as the limiting processes in testing for normality with estimated mean and/or variance. 

Such processes require concrete fine analysis. For the Kac--Kiefer--Wolf\-owitz processes, the exact $L_2$-small ball asymptotics was calculated in~\cite{NazPet2015-SVFpert}, for several other important Durbin-type processes it was managed in~\cite{Pet2020}. The main tools in obtaining these results are the asymptotic expansion of oscillatory integrals with amplitude being a slowly varying function~\cite[Theorems 1--3]{NazPet2015-SVFpert} and the following generalization of the formula \eqref{eq:P-model-naznik}:
\smallskip

{\it
Suppose that $\th>0$, $\delta>-1$, $p>1$. 
Let a function $\Phi$ be slowly varying at infinity, monotonically tending to zero, and assume that the integral $\int\limits_1^{\infty}\frac{\Phi(t)}{t}\,dt$ diverges.\footnote{If this integral converges then the asymptotics up to a constant coincides with \eqref{eq:P-model-naznik}.}
Then, as $\ep\to0$,
\begin{multline*}
\mathbb{P}\Big\{ \sum\limits_{k=1}^{\infty}(\th (k+\delta+\Phi(k)))^{-p}\xi_k^2\leq\ep^2  \Big\}\\
\sim C\cdot \ep^{\gamma} \exp\bigg( -D \ep^{-\frac{2}{p-1}} +\frac p2 \int\limits_1^{\ \ep^{-\frac{2}{p-1}}}\frac{\Phi(t)}{t}\,dt\bigg),
\end{multline*}
where $D$ and $\gamma$ are defined in~\eqref{eq:D-gamma} while $C=C(\th,\delta,d,\Phi)$ is some (unknown) constant.
}

\subsubsection{Detrended Green Gaussian processes}

This class of processes generalizes the notion of the demeaned process. It is natural to view the process $\overline{X}(t)$ in~\eqref{eq:integrated-centered-process} as the projection of $X(t)$ onto the subspace of functions orthogonal to the constants in $L_2[0,1]$. In particular, this means
\begin{align*}
\| \overline{X}\|_{2,[0,1]}=\min\limits_{a\in\mathbb{R}} \|X-a\|_{2,[0,1]}.
\end{align*}

In a similar way, we can define the $n$-th order \textit{\textbf{detrended}} process $X_{\langle n\rangle}(t)$ as the projection of $X(t)$ onto the subspace of functions orthogonal in $L_2[0,1]$ to the set ${\cal P}_n$ of polynomials of degree up to $n$, that is
\begin{align*}
\| X_{\langle n\rangle}\|_{2,[0,1]}=\min\limits_{P\in{\cal P}_n} \|X-P\|_{2,[0,1]}.
\end{align*}
This process is given by the formula
\begin{align*}
X_{\langle n\rangle}(t)=X(t)-\sum\limits_{j=0}^n a_jt^j,
\end{align*}
where the random variables $a_j$ are determined by relations
\begin{align*}
\int\limits_0^1 t^jX_{\langle n\rangle}(t)\,dt=0,\quad j=0,\ldots,n.
\end{align*}
The first order detrended processes are simply called detrended.
\medskip

For several concrete demeaned Green  Gaussian processes on $[0,1]$, the $L_2$-small ball asymptotics were obtained in~\cite{BegNikOrs2005} and~\cite{Ai2016}\footnote{Notice that the power term in~\cite[Proposition~2]{Ai2016} was calculated erroneously.}; some particular detrended processes were managed in~\cite{AiLiLiu2012} (up to a constant) and~\cite{KirNik2014}.

Yu.P.~Petrova \cite{Pet2017} considered this problem in a more general setting. Let $X$ be a Green Gaussian process on $[0,1]$ corresponding to the boundary value problem~\eqref{eq:L-bvp} with ${\cal L}=(-1)^\ell u^{(2\ell)}$. It turns out that for $n\ge 2\ell$ the covariance function of the process $X_{\langle n\rangle}$ does not depend on original boundary conditions\footnote{Notice that the notation $X_{\langle n\rangle}$ in~\cite{Pet2017} corresponds to $X_{\langle n-1\rangle}$ here.} and is the generalized Green function of the problem
\begin{align*}
(-1)^\ell u^{(2\ell)}(t)
=
\lambda u(t)+P(t),
&& 	
\int\limits_0^1 t^j u(t)\, dt=0,\quad j=0,\dots, n,
\end{align*}
where $P\in{\cal P}_{n-2\ell}$ is a polynomial with unknown coefficients.

The final result in~\cite{Pet2017} reads as follows:\footnote{For $X=W$ and $X=B$, the eigenvalues of $G_{X_{\langle n\rangle}}$ were found earlier in~\cite{AiLi2014}. However, the power term in $L_2$-small ball asymptotics was calculated erroneously in this paper.}
\smallskip

{\it
Let $n\ge 2\ell$. Then, as $\ep\to0$,
\begin{align*}
\mathbb{P}\Big\{\|X_{\langle n\rangle}\|_2\leq\ep\Big\}
\sim 
C
\ep^{\gamma}
\exp\Bigl(-D\ep^{-\frac2{2\ell-1}}\Bigr), 
\end{align*}
where $D$ is defined in~\eqref{eq:gamma} with $\theta_\ell=1$, $\gamma=\frac{(\ell-1)^2-2n\ell}{2\ell-1}$, while $C=C(\ell, n)$ is an explicit constant.
}

\subsection{Fractional Gaussian noise and related processes}
It is well known that the standard Brownian motion (the Wiener process) is a primitive of the Gaussian white noise, that is the zero mean distribution-valued Gaussian process with the identity covariance operator.

Now we recall that the \textbf{\textit{fractional Brownian motion (fBm)}} $W^H$ with the \textbf{\textit{Hurst index}} $H\in(0,1)$ is the zero mean Gaussian process with the covariance function
\begin{align}
\label{eq:frac-BM}
G_{W^{H}}(t,s)=\frac12\,\bigl(s^{2H}+t^{2H}-|s-t|^{2H}\bigr),\qquad t,s>0
\end{align}
(for $H=\frac 12$ we obtain the conventional Wiener process).

In a similar way, the fractional Brownian motion (say, on the interval $[0,1]$) can be considered as a primitive of the \textbf{\textit{fractional Gaussian noise}} (see, e.g.,~\cite{Qian2003}) 
that is the zero mean distribution-valued Gaussian process with the covariance operator ${\cal K}_{2-2H}$, where
\begin{align*}
({\cal K}_{\alpha}u)(x)=(1-\tfrac\alpha2)\, \frac d{dx}\int\limits_0^1 \text{sign}(x-y)|x-y|^{1-\alpha}u(y)\,dy.
\end{align*}
It is easy to check that ${\cal K}_1=\mathrm{Id}$, as required.
\medskip

The logarithmic $L_2$-small ball asymptotics for the fBm and similar processes was extensively studied at the beginning of XXI century (see Section~\ref{sec:logarithm}). In contrast, corresponding exact $L_2$-small ball asymptotics was a challenge until recently. This is related to the fact that $W^H$ is not a Green Gaussian process for $H\ne\frac 12$.

The breakthrough step in this problem was managed by P.~Chigansky and M.~Kleptsyna~\cite{ChiKlep2018}. Using the substitution $u(t)=\int\limits_t^1\varphi(\tau)\,d\tau$, the equation 
\begin{align*}
\int\limits_0^1 G_{W^{H}}(t,s)\varphi(s)\,ds=\mu\varphi(t), \qquad 0\le t\le 1,
\end{align*}
for eigenvalues of the covariance operator of fBm 
was reduced to the \textbf{\textit{generalized eigenproblem}} for the second order ODO (here $\lambda=\mu^{-1}$)
\begin{equation}
\label{eq:gen-eigenpr}
-u''=\lambda {\cal {K}}_{\alpha}u \quad {\rm
on}\quad [0,1];\qquad u'(0)=u(1)=0.
\end{equation}
By the Laplace transform the problem (\ref{eq:gen-eigenpr}) was converted to the Riemann--Hilbert problem which, in turn, was solved asymptotically using the idea of S.~Ukai~\cite{Ukai1971} and B.V.~Pal'tsev \cite{Pal1974}, \cite{Pal2003}. In this way the two-term asymptotics of the eigenvalues with the remainder estimate was obtained for the fBm in the full range of the Hurst index. Based on this, the exact $L_2$-small ball asymptotics 
for $W^H$ was established for the first time, along with some other applications. It should be mentioned that the eigenfunctions asymptotics for fBm was also obtained in \cite{ChiKlep2018}. 

In later papers \cite{ChiKlepMar2020}, \cite{ChiKlepMar2020-bridges}, \cite{KlepMarChi2020} similar results were obtained for some other particular fractional Gaussian processes. 

A.I. Nazarov~\cite{Naz2021-frac} considered the problem in a more general setting. He considered the generalized eigenproblem 
\begin{equation}
\label{eq:gen-eigenpr-1}
-u''+p_0u=\lambda {\cal {K}}_{\alpha}u \quad {\rm
on}\quad [0,1]
\end{equation}
(here $p_0\in L_1[0,1]$)
with general self-adjoint boundary conditions. 

It turns out that the method invented in~\cite{ChiKlep2018} runs without essential changes in general case. Moreover, the term $p_0u$ in (\ref{eq:gen-eigenpr-1}) can be considered as a weak perturbation which does not affect the two-term eigenvalues asymptotics.\footnote{The latter fact follows from the general result \cite{Naz2020-lemmata} on spectrum perturbations and careful analysis of the eigenfunctions asymptotics based on ideas of~\cite{ChiKlep2018}.} 
\medskip

The final result in~\cite{Naz2021-frac} reads as follows:
\smallskip

{\it Let the equation for eigenvalues of the covariance operator of the Gaussian process $X$ on $[0,1]$ can be reduced to the generalized eigenproblem~\eqref{eq:gen-eigenpr-1}. If the boundary conditions do not contain the spectral parameter\footnote{We do not know examples of conventional Green Gaussian processes corresponding to a boundary value problem with the spectral parameter in the boundary conditions. In contrast, some natural fractional Gaussian processes arrive at the problem \eqref{eq:gen-eigenpr-1} with $\lambda$ in boundary conditions, see \cite[Section 4]{Naz2021-frac}. In this case the exponent $\gamma$ in~\eqref{eq:P-fractional-green} differs from~\eqref{eq:DH-gamma} and should be derived separately though the general method~\cite{Naz2021-frac}  still works.} then, as $\ep\to0$,
\begin{align}
\label{eq:P-fractional-green}
\mathbb{P}\{{}\|X\|_{2,[0,1]}\le \ep\}\sim {\cal C}(X)\cdot 
\ep^{\gamma} \exp\Big(-D(H)\ep^{-\frac {1}{H}} \Big).
\end{align}
where
\begin{align}
\label{eq:DH-gamma}
D(H)=\tfrac H{(2H+1)\sin\left(\frac {\pi}{2H+1}\right)} \Big(\tfrac {\sin(\pi H)\Gamma(2H+1)} {(2H+1)\sin\left(\frac {\pi}{2H+1}\right)}\Big)^{\frac {1}{2H}},
\qquad
\gamma=\frac{(H-\frac12)^2+\varkappa}{2H},
\end{align}
(recall that $\varkappa$ is the sum of orders of normalized boundary conditions in~\eqref{eq:gen-eigenpr-1})
and ${\cal C}(X)$ is some (unknown) constant.
 If the problem~\eqref{eq:gen-eigenpr-1} has a zero eigenvalue then formula~\eqref{eq:P-fractional-green} should be corrected similarly to Remark~\ref{rmk:zero-eigen-1}. 
}
\medskip

The result of \cite{Naz2021-frac} encompasses some previous results and covers several more fractional Gaussian processes.\footnote{Notice that there exist various fractional analogs of classical Gaussian processes (Brownian bridge, Ornstein--Uhlenbeck process etc.). Only for some of them the $L_2$-small ball asymptotics up to a constant can be managed by the considered method. In contrast, the algorithm of derivation of logarithmic asymptotics considered in Section~\ref{sec:logarithm} covers all this variety.}

It is worth to note the paper  \cite{ChiKlep2021}, where the two-term eigenvalues asymptotics for the Riemann--Liouville process and the Riemann--Liouville bridge on $[0,1]$ was derived by reduction to the generalized eigenproblem for \textit{\textbf{fractional differential operators}} 
\begin{align}
\label{eq:RL-Cap}
D^\alpha_{1-}D^\alpha_{0+}u=\lambda u\qquad \mbox{and} \qquad {^c}D^\alpha_{1-} {^c}D^\alpha_{0+}u=\lambda u
\end{align}
with proper boundary conditions (here $D^\alpha_{0+}$,  $D^\alpha_{1-}$, and ${^c}D^\alpha_{0+}$, ${^c}D^\alpha_{1-}$ are the left/right Riemann--Liouville and Caputo fractional derivatives, respectively, see \cite[Chapter 2]{KilSriTru2006-book}).
\medskip

The next natural task is to transfer the result of~\cite{ChiKlep2018} to the (multiply) integrated fractional processes. The first step here was made by 
P.~Chigansky, M.~Kleptsyna and D.~Marushkevych~\cite{ChiKlepMar2018} who derived the two-term eigenvalues asymptotics and the $L_2$-small ball asymptotics up to a constant for the integrated fBm $(W^H)_1$ on $[0,1]$.

A unified approach similar to~\cite{Naz2021-frac} reduces the problem for $m$-times integrated fBm and related processes to the generalized eigenproblem
\begin{align}
\label{eq:gen-eigenpr-2}
    {\cal L}u=\lambda {\cal {K}}_{\alpha}u \quad {\rm
on}\quad [0,1], \qquad u\in {\cal D}({\cal L}),
\end{align}
where ${\cal L}$ is an ODO of the form~\eqref{eq:L} with $\ell=m+1$ and $p_\ell\equiv1$, whereas the domain ${\cal D}({\cal L})$ is defined by the boundary conditions \eqref{eq:bc-non-separated} of special form. More generally, if the equation for eigenvalues of the covariance operator of the Gaussian process $X$ on $[0,1]$ can be reduced to the generalized eigenproblem~\eqref{eq:gen-eigenpr-2}, we call $X$
the \textit{\textbf{fractional-Green Gaussian process}}.

\medskip
\textbf{Conjecture.} 
\textit{Let $X$ be a fractional-Green Gaussian process on $[0,1]$. If the boundary conditions do not contain the spectral parameter then, as $\ep\to0$,
\begin{align}
\label{eq:P-gen-fractional-green}
\mathbb{P}\{{}\|X\|_{2,[0,1]}\le \ep\}\sim {\cal C}(X)\cdot 
\ep^{\gamma} \exp\Big(-D_\ell(H)\ep^{-\frac {1}{\ell-1+H}} \Big),
\end{align}
where
\begin{align}
\label{eq:Dl-fractional-green}
D_\ell(H)&=\tfrac
{\ell-1+H}{(2(\ell+H)-1)
\sin\left(\frac
{\pi}{2(\ell+H)-1}\right)}
\cdot 
\left(\tfrac
{\Gamma(2H+1)\sin(\pi H)}{(2(\ell+H)-1) \sin\left (\frac
{\pi}{2(\ell+H)-1}\right)}\right)^{\frac {1}{2(\ell-1+H)}},
\end{align}
\begin{align}
\label{eq:gamma-fractional-green}
\gamma&=-\frac{(\ell-1)(2\ell+1)}{2(H+\ell-1)}+\frac{(H-\frac12)^2+\varkappa}{2(H+\ell-1)},
\end{align}
$\varkappa$ is the sum of orders of normalized boundary conditions in~\eqref{eq:gen-eigenpr-2},
and ${\cal C}(X)$ is some (unknown) constant. If the problem~\eqref{eq:gen-eigenpr-1} has zero eigenvalues then formula~\eqref{eq:P-fractional-green} should be corrected similarly to Remark~\ref{rmk:zero-eigen-1}. 
}
\medskip

At the moment, this conjecture is proved for operators ${\cal L}$ with separated boundary conditions \cite{ChiKlepNazRast2023}. Notice that the problem of derivation of the exact constant ${\cal C}_X$ in $L_2$-small ball asymptotics for the fractional Gaussian processes is completely open.
\medskip

Formulae~\eqref{eq:P-gen-fractional-green}--\eqref{eq:gamma-fractional-green} imply that, as for the Green Gaussian processes (see Corollary~\ref{cor:integrated-constant}),
for any given fractional-Green Gaussian process $X$ on $[0,1]$, the $L_2$-small ball asymptotics for various integrated processes $X_m^{[\beta_1,\ldots,\beta_m]}$ can differ only by a constant. To find a more general class of Gaussian processes with such property is an interesting open problem.

\subsection{Tensor products}
\label{subsec:tensor-exact}
In this subsection we consider multiparameter Gaussian functions (Gaussian random fields) of ``tensor product'' type. It means that the covariance function of this random field can be decomposed in a product of ``marginal'' covariances depending on different arguments. The classical examples of such fields are the Brownian sheet, the Brownian pillow and the Brownian pillow-slip (the Kiefer field), see, e.g.,~\cite{VaarWell1986-book}. Less known examples can be found in~\cite{Kho2002-book} and~\cite{LifNazNik2003}. The notion of tensor products of Gaussian processes or Gaussian measures was known long ago, see~\cite{Car1977} and~\cite{CarChe1979} for a more general approach. Such Gaussian fields are also studied in related domains of mathematics, see, e.g., \cite{PapWas1990}, \cite{RitWasWoz1995}, \cite{GraLusPag2003}, \cite{LusPag2004} and \cite{KhaLim2022}.
\medskip

Suppose we have two Gaussian random functions $X(x)$, $x \in \mathbb{R}^m$, and $Y(y)$, $y \in \mathbb{R}^n$, with zero means
and covariance functions $G_X(x,u)$, $x,u \in \mathbb{R}^m$, and
$G_Y(y,v)$, $y,v \in \mathbb{R}^n$, respectively. Consider
the new Gaussian function $Z(x,y)$, $x \in \mathbb{R}^m$, $y \in \mathbb{R}^n$, which has zero mean and the covariance function
$$
G_Z((x,y),(u,v)) =G_X(x,u)G_Y(y,v).
$$
Such Gaussian function obviously exists, and corresponding covariance operator $\mathcal G_Z$ is the tensor product of two ``marginal'' covariance operators, that is $\mathcal G_Z=\mathcal G_X\otimes \mathcal G_Y$. So we use in the sequel the notation $Z = X\otimes Y$ and 
we call the Gaussian field $Z$ \textit{\textbf{the tensor product}} of the fields $X$ and $Y$ using sometimes for shortness the term ``field-product''. The generalization to the multivariate case when obtaining the fields $\underset{j=1}{\overset{d}\otimes} X_j$ 
is straightforward. It is easy to see that the tensor products 
$$
\mathbb{W}^d(x_1,\dots,x_d)=\underset{j=1}{\overset{d}\otimes} W(x_j)\qquad\text{and}\qquad \mathbb{B}^d(x_1,\dots,x_d)=\underset{j=1}{\overset{d}\otimes} B(x_j),
$$ 
are, respectively, just the $d$-dimensional Brownian sheet and the $d$-di\-men\-sion\-al Brownian pillow.

\begin{remark}
\label{rm:tensor-sums}
We underline that the eigenvalues of the operator ${\cal K}_1\otimes {\cal K}_2$ are just pairwise products of the eigenvalues of the marginal operators ${\cal K}_1$~and~${\cal K}_2$. Therefore, formula \eqref{eq:KL} for the fields-products deals in fact with multiple sums.
\end{remark}

The logarithmic small ball asymptotics for fields-products are widely investigated (see Subsection~\ref{subsec:reg-var}), whereas examples on exact $L_2$-small ball asymptotics for fields-products are not numerous. The first results were obtained by J.~Fill and F.~Torcaso~\cite{FillTor2004}, 
who considered multiply integrated Brownian sheets
$$
\mathbb{W}^d_{m_1,\dots,m_d}(x_1,\dots,x_d)=\underset{j=1}{\overset{d}\otimes} W_{m_j}^{[0,\dots,0]}(x_j)
$$
(see formula~\eqref{eq:m-times-interated-X}) on the unit cube ${\cal Q}_d=[0,1]^d$. Using the Mellin transform, authors succeeded to obtain the \textit{full asymptotic expansion} of the logarithm of the small ball probability as $\ep\to0$.

It turns out that to extract from this expansion the exact small ball asymptotics, say, for the conventional Brownian sheet $\mathbb{W}^d$, one needs infinite number of terms. However, in some cases finite number is sufficient. We restrict ourselves to only one such example:\footnote{Notice that in~\cite[Example~5.9]{FillTor2004} the minus sign is lost under the exponent.}
$$
\mathbb{P}\{\|\mathbb{W}^2_{0,m}\|_{2,{\cal Q}_2}\le \ep\}\sim C\ep\exp\big(-D_1\ep^{-2}-D_2\ep^{-\frac 2{m+1}}\big),\quad \ep\to0,
$$
where $C$, $D_1$, $D_2$ are explicit (though complicated) constants.
\medskip

Another technique based on direct inversion of the Laplace transform was used by L.V.~Rozovsky~\cite{Roz2019-fields}, \cite{Roz2019}. He has obtained the exact $L_2$-small ball asymptotics for double sums (cf. Remark~\ref{rm:tensor-sums})
\begin{align}
\label{eq:double-sum}
{\cal S}^{(2)}:=a^2\sum_{j=1}^{\infty}\sum_{k=1}^{\infty}(j+\delta_1)^{-p_1}(k+\delta_2)^{-p_2}\xi_{jk},
\end{align}
where $p_1>1$, $p_2=2$, $\delta_1,\delta_2>-1$, whereas $\xi_{jk}$ are i.i.d. standard normal random variables.\footnote{In~\cite{Roz2019-fields} the author deals with the case $p_1\ne2$, in~\cite{Roz2019} --- with the case $p_1=2$.}

\section{Logarithmic asymptotics}
\label{sec:logarithm}

\subsection{General assertions}
The first general result on logarithmic $L_2$-small ball asymptotics was obtained by W.V.~Li~\cite{Li1992} who proved that for arbitrary $N\in\mathbb{N}$ we have, as $\ep\to 0$,
\begin{align*}
\log \mathbb{P}\Big\{\sum_{k=1}^{\infty} \mu_k \xi_k^2\leq\ep^2\Big\}
\sim 
\log \mathbb{P}\Big\{\sum_{k=N}^{\infty} \mu_k \xi_k^2\leq\ep^2\Big\}.
\end{align*}
This relation, in particular, implies that the logarithmic asymptotics does not change under any finite-dimensional perturbation of a Gaussian random function.\footnote{This follows from the fact that by the minimax principle (see, e.g., \cite[Section~9.2]{BirSol2010-book}) the eigenvalues $\widetilde \mu_k$ of the perturbed operator satisfy the two-sided estimate 
$$
\mu_{k+N}\le \widetilde \mu_k \le \mu_{k-N},
$$ 
where $\mu_k$ are the eigenvalues of the original (compact, self-adjoint) operator, and $N$ is the rank of perturbation.}

An important development of this result is the so-called \textit{\textbf{log-level comparison principle}}. Roughly speaking, it means that if the eigenvalues of covariance operators for two Gaussian vectors have the same one-term asymptotics then the logarithmic $L_2$-small ball asymptotics for these vectors coincide. In particular, the power (regularly varying) eigenvalues asymptotics yields the power (respectively, regularly varying) logarithmic $L_2$-small ball asymptotics.

For concrete behavior of eigenvalues (power decay, regular decay, etc.) such statements were proved in various papers, see, e.g., \cite{NazNik2005-frac}, \cite{KarNazNik2008-Tensor}, \cite{Naz2004-logselfsim}. In general case the result was established independently by F.~Gao and W.V.~Li~\cite{GaoLi2007-log} and A.I.~Nazarov~\cite{Naz2009-log}.\footnote{The statement in~\cite{GaoLi2007-log} is more general than in~\cite{Naz2009-log} but the assumptions on the eigenvalues behavior are somewhat more restrictive.}
\medskip

We recall the notion of \textit{\textbf{the counting function}} of the sequence $(\mu_k)$:
\begin{align*}
{\cal N}(\tau):=\#\{k:\mu_k>\tau\}.
\end{align*}
This notion is common in spectral theory. Notice that the functions  $k\mapsto \mu_k$ and $\tau\mapsto {\cal N}(\tau)$ are in essence mutually inverse. So, if $\mu_k$ decay not faster than a power of $k$ and have ``not too wild'' behavior\footnote{\label{foot:count-eigen}For instance, this is true for regular behavior of $\mu_k$. Moreover, a sequence $(\mu_k)$ is regularly varying with index $\rho<0$ if and only if its counting function is regularly varying at zero with index $\frac 1\rho$.} then the
one-term asymptotics of ${\cal N}(\tau)$ as $\tau\to0$ provides the one-term asymptotics of $\mu_k$ as $k\to\infty$, and vice versa.\footnote{The latter is not true if $\mu_k$ decay slower than any power of $k$, but in our case it is impossible since $\mu_k$ are summable.} However, for the super-power decreasing of $\mu_k$ to obtain the asymptotics of the counting function is much simpler than the asymptotics of eigenvalues.\medskip

Now we are ready to formulate the log-level comparison principle in the form of~\cite{Naz2009-log}. 

\begin{theorem}[{\cite[Theorem 1]{Naz2009-log}}]
\label{thm:log-comparison}	
Let $(\mu_k)$ and
$(\widetilde\mu_k)$ be positive summable
sequences with counting functions ${\cal N}(\tau)$ and $\widetilde{\cal N}(\tau)$, respectively. Suppose that the function ${\cal N}$ satisfies
\begin{equation}
\label{eq:not-too-slow}
\liminf\limits_{x\to0}\frac {\int\limits_0^{hx} {\cal N}(\tau)\,d\tau}
{\int\limits_0^{x}{\cal N}(\tau)\,d\tau}>1 \qquad
\mbox{for any} \qquad h>1. 
\end{equation}
If
$$
{\cal N}(\tau)\sim\widetilde{\cal N}(\tau), \qquad \tau\to 0,
$$
then, as $\ep\to0$,
\begin{equation}
\label{eq:log-level}
\log\mathbb{P}\Big\{\sum_{k=1}^\infty\widetilde\mu_k\xi_k^2
\le \ep^2 \Big\} \sim \log\mathbb{P}\Big\{ \sum_{k=1}^\infty
\mu_k\xi_k^2\le \ep^2 \Big\}.
\end{equation}
\end{theorem}

\begin{remark}
The assumption \eqref{eq:not-too-slow} is satisfied, for instance, for ${\cal N}(\tau)$ regularly varying at zero with index $\rho\in(-1,0]$. However, if ${\cal N}(\tau)= C\tau^{-1}\log^{-\sigma}(\tau^{-1})$, $\sigma>1$, then the relation
\eqref{eq:log-level} in general fails as it is shown in \cite[Section~4]{KarNazNik2008-Tensor}.\footnote{For this case, the \textit{\textbf{double logarithmic}} asymptotics was derived in~\cite[Proposition~4.4]{KarNazNik2008-Tensor}:
$$ 
\log \Bigl|\log\mathbb{P}\Big\{ \sum_{k=1}^\infty
\mu_k\xi_k^2\le \ep^2 \Big\}\Bigr| \sim
 \left (\frac {C}{\sigma-1}\right)^{\frac 1{\sigma-1}}\ep^{-\frac
{2}{\sigma-1}}, \qquad \ep \to 0. 
$$
}  
Thus, the assumption~\eqref{eq:not-too-slow}  cannot be removed.
\end{remark}

Theorem~\ref{thm:log-comparison} allows to use numerous results in spectral theory of integral operators, where the counting functions for eigenvalues and their asymptotics are widely investigated. In the following subsections we provide the results related to various behavior of ${\cal N}(\tau)$ as $\tau\to0$.

\begin{remark}
It is of great importance that the one-term eigenvalues asymptotics is quite stable under perturbations. We provide here two assertions for the additive perturbations of a compact operator or of its inverse.

For the operators with power eigenvalues asymptotics the first assertion was established by H.~Weyl~\cite{Weyl1912}; see also~\cite[Lemmata 1.16 and 1.17]{BirSol1974-10school}. 
    
\begin{theorem}
\label{thm:spectral-stability}
Let ${\cal K}$ be a compact self-adjoint positive operator, and let the counting function of its eigenvalues ${\cal N}_{\cal K}(\tau)$ have the following asymptotics:
$$
{\cal N}_{\cal K}(\tau)\sim \Phi(\tau),\quad \tau\to0,
$$
where the function $\Phi$ is regularly varying with index $\rho<0$ at zero.

{\bf 1} (\cite[Lemma 3.1]{BirSol2001}). Assume that ${\cal K}_1$ is another compact self-adjoint positive operator, and that ${\cal N}_{{\cal K}_1}(\tau)=o(\Phi(\tau))$ as $\tau\to0$. Then
$$
{\cal N}_{{\cal K}+{\cal K}_1}(\tau)\sim {\cal N}_{{\cal K}-{\cal K}_1}(\tau)\sim {\cal N}_{\cal K}(\tau),\quad \tau\to0.
$$

{\bf 2} (\cite[Lemma 5.1]{NazShei2012}).\footnote{In~\cite{NazShei2012} the assumptions on $\Phi$ are somewhat weaker. Notice that a much more general statement was proved in~\cite[Theorem 3.2]{MarMat1982}.} Assume that ${\cal B}$ is a self-adjoint operator such that ${\cal K}{\cal B}$ is compact, and ${\cal K}^{-1}+{\cal B}$ is positive definite. Then
$$
{\cal N}_{({\cal K}^{-1}+{\cal B})^{-1}}(\tau)\sim {\cal N}_{\cal K}(\tau),\quad \tau\to0.
$$
\end{theorem}

\end{remark}

\subsection{Power asymptotics}

Some explicit formulae of power logarithmic asymptotics for concrete (Green) Gaussian processes were derived in~\cite{Li1992}, \cite{KhoShi1998} and \cite{ChenLi2003}.\footnote{Later these results were covered by the exact $L_2$-small ball asymptotics of the Green Gaussian processes, see Theorem~\ref{thm:small-ball-Green}.}

The first result on logarithmic asymptotics for a non-Green Gaussian process (namely, for the fractional Brownian motion, see~\eqref{eq:frac-BM}) was obtained by 
J.~Bronski~\cite{Bro2003}. By particular fine analysis, he derived the one-term eigenvalues asymptotics with the remainder estimate, and therefore, the relation
\begin{align}
\label{eq:P-log-fBm}
\log\mathbb{P}\{{}\|W^H\|_{2,[0,1]}\le \ep\}\sim -\,D(H)\ep^{-\frac {1}{H}}, \qquad \ep\to0,
\end{align}
where $D(H)$ is given in~\eqref{eq:DH-gamma}.
\medskip

A.I. Nazarov and Ya.Yu. Nikitin~\cite{NazNik2005-frac} 
succeeded to manage quite general class of Gaussian random functions (including weighted, multiparameter and vector-valued ones). They used a particular case of the deep results by M.S.~Birman and M.Z.~Solomyak~\cite{BirSol1970} (see also~\cite[Appendix 7]{BirSol1974-10school}) on spectral asymptotics of integral operators. It turns out that if the covariance function has the main singularity of power type on the diagonal then the eigenvalues $\mu_k$ have the power one-term asymptotics. Also the argument in~\cite{NazNik2005-frac} uses part {\bf 1} of Theorem~\ref{thm:spectral-stability}, a particular case of Theorem~\ref{thm:log-comparison}, and the following corollary of the formula~\eqref{eq:P-model-naznik}: for $\th>0$ and $p>1$,
\begin{align}
\label{eq:power-asymp-small-balls}
\log\mathbb{P}\Big\{\sum_{k=1}^{\infty}(\th k)^{-p} \xi_k^2 \leq
\ep^2\Big\}
\sim -\, D \ep^{-\frac {2}{p-1}},\qquad \ep\to0,
\end{align}
where $D$ is given in~\eqref{eq:D-gamma}.

The final result of~\cite{NazNik2005-frac} reads as follows.\footnote{Here it is given in a simplified form. In particular, we restrict ourselves to the one-parameter and scalar-valued case.}

\begin{theorem}
\label{thm:log-asymp-singularity-power}
Let $X$ be a zero mean Gaussian process on $[0,1]$, and let the covariance function admit the representation
\begin{align*}
G_X(t,s)=f(|t-s|^{2\rho})\,{\cal R}(t,s) + \sum\limits_{j=1}^N a_j(t)a_j(s),\qquad t,s\in [0,1].
\end{align*}
Assume that $\rho>0$, $\rho \notin{\mathbb
N}$, $a_j\in L_2[0,1]$, ${\cal
R}(t,s)\equiv{\cal R}(s,t)$ is a 
smooth function on $[0,1]\times [0,1]$
 while $f$ is a smooth function on $\mathbb{R}_+$ with $|f'(0)|=1$.\footnote{The latter assumption can be made without loss of generality. We also notice that $R(t,t)$ has constant sign for $t\in [0,1]$. }

 Then for any non-negative function $\psi\in L_2[0,1]$ we have, as $\ep\to0$,
\begin{align}
\label{P-log-naznik}
\log\mathbb{P}\{{}\|\psi X\|_{2,[0,1]}\le \ep\}\sim - \,\tfrac {\rho\,\mathfrak{J}_\rho}{(2\rho+1)\sin\left(\frac {\pi}{2\rho+1}\right)} \Big(\tfrac {2\,\mathfrak{J}_\rho\,\Gamma(2\rho+1)|\sin(\pi \rho)|} {(2\rho+1)\sin\left(\frac {\pi}{2\rho+1}\right)}\Big)^{\frac {1}{2\rho}}\ep^{-\frac {1}{\rho}},
\end{align}
where 
$$
\mathfrak{J}_\rho=\int\limits_0^1\big(|{\cal R}(t,t)|\,\psi^2(t)\big)^{\frac1{2\rho+1}}\,dt.
$$
\end{theorem}

For $X=W^H$ we have $\rho=H$, ${\cal R}(t,s)\equiv-\frac 12$, and for $\psi(t)\equiv 1$ formula \eqref{P-log-naznik} coincides with \eqref{eq:P-log-fBm}.

\begin{remark}
    Theorem~\ref{thm:log-asymp-singularity-power}, in particular, implies that the terms more smooth than the main (the most singular) one do not influence upon the logarithmic $L_2$-small ball asymptotics.\footnote{This remark shows that the result in~\cite{CaiHuWang2021} is incorrect.}
\end{remark}

If our Gaussian process is stationary, it is natural to derive its small ball asymptotics in terms of its spectral measure. Such results were obtained by M.A.~Lifshits and A.I.~Nazarov \cite{LifNaz2018}. Authors considered weighted stationary sequences (for this case, see also the earlier paper~\cite{HongNazLif2016} with S.Y.~Hong) and weighted stationary processes on the line and on the circle. Along with real processes, the class of proper complex processes was treated.

We provide one of the results of~\cite{LifNaz2018}.

\begin{theorem}[{\cite[Theorem~2.2]{LifNaz2018}}]    
Let $X(t)$ be a (real) $2\pi$-periodic zero mean,
mean-square continuous stationary Gaussian process.
Assume that its spectral measure\footnote{Under our assumptions it is a symmetric non-negative sequence $\mu_k=\mu_{-k}$, $k\in\mathbb{Z}$.}  satisfies the asymptotic condition
\begin{align*}
\mu_k \sim M |k|^{-r}, \qquad  |k|\to \infty,
 \end{align*}
with some $r>1$ and $M>0$. Let $\psi\in L_2[0,2\pi]$.

Then we have, as $\ep\to 0$,
\begin{align*}
    \log \mathbb{P}\left( \|\psi X\|_{2,[0,2\pi]} \le \ep\right)
\sim -\, \frac{r-1}{2} \bigg(\frac{(2\pi M)^{\frac 1r}}{r\sin(\frac\pi{r})}
\int\limits_{0}^{2\pi} \psi^{\frac 2{r}}(t)\,  dt\bigg)^{\frac{r}{r-1}}
 \ep^{-\frac{2}{r-1}}\, .
 \end{align*}
\end{theorem}

\noindent The results of~\cite{HongNazLif2016} and~\cite{LifNaz2018} are based on the spectral theory of pseudo-differential operators with homogeneous symbols developed by M.S.~Birman and M.Z. Solomyak~\cite{BirSol1977}, \cite{BirSol1979}.

To generalize the results of~\cite{LifNaz2018} for stationary processes on general Abelian groups is an interesting open problem.

\subsection{Regularly varying asymptotics}
\label{subsec:reg-var}

We begin with the generalization of the relation~\eqref{eq:power-asymp-small-balls} for the regularly varying asymptotics:

\begin{lemma}[{\cite[Theorem~4.2 and Proposition~4.3]{KarNazNik2008-Tensor}}] 
\label{lm:log-P-reg-spect-asymp}
Let $p>1$. Suppose that the function $\Psi(t)$ is slowly varying at infinity, and\,\footnote{The relations~\eqref{eq:equiv-slow} can be assumed without loss of generality, since for any SVF $\Phi$ there is a SVF $\Psi$ satisfying~\eqref{eq:equiv-slow} such that $\Psi(t)\sim \Phi(t)$ as $t\to+\infty$, see, e.g.,~\cite[Chapter~1]{Sen1976-book}.}
\begin{equation}
\label{eq:equiv-slow}
t\cdot(\ln(\Psi))'(t)\to 0,\qquad
t^2\cdot(\ln(\Psi))''(t)\to 0,\qquad t\to+\infty.
\end{equation}
Then, as $\ep\to0$,
\begin{align}
\label{eq:reg-asymp-small-balls}
\log\mathbb{P}\Big\{\sum_{k=1}^{\infty}k^{-p}\Psi(k) \,\xi_k^2 \leq
\ep^2\Big\}
\sim - \,\ep^{-\frac {2}{p-1}}\Theta(\ep^{-1}),
\end{align}
where $\Theta$ is an explicit (though complicated) SVF depending on $p$ and on $\Psi$. In particular, if $\th >0$ and $\sigma\in\mathbb{R}$ then
\begin{align*}
\log\mathbb{P}\Big\{\sum_{k=1}^{\infty}(\th k)^{-p} \log^{\sigma}(k+1)\,\xi_k^2 \leq
\ep^2\Big\}
\sim - \,D\ep^{-\frac {2}{p-1}}\log^{\frac{\sigma}{p-1}}(\ep^{-1}), \quad \ep\to0,
\end{align*}
where $D=\big(\frac {p-1}{2}\big)^{1-\frac{\sigma}{p-1}}\Big( \frac
{\pi}{\th p\sin\frac{\pi}{p}}\Big)^{\frac {p}{p-1}}$.
\end{lemma}

A simple example of Gaussian process with regular eigenvalues behavior (and therefore regularly varying $L_2$-small ball asymptotics) can be constructed as follows (see, e.g., \cite{KarNazNik2008-Tensor}). Let $X(t)$ be the stationary Gaussian process on $[0,1]$ with zero mean and covariance function
$$
G_X(t,s)=g(|t-s|),
$$ 
where $g$ is a smooth function on $\mathbb{R}_+$, and 
\begin{align*}
g(t)-g(0)=-t^\alpha\Psi(t^{-1})+O(t^\beta),\qquad t\to0+,
\end{align*}
where $0<\alpha<2$, $\beta>\alpha$, while $\Psi$ is a SVF. Then, similarly to \cite[Appendix~7]{BirSol1974-10school}, one can show that
\begin{align*}
    \mu_k\sim \frac{2\sin\bigl(\frac{\pi\alpha}{2}\bigr)\Gamma(\alpha+1)}{\pi^{\alpha+1}}\cdot\frac{\Psi(k)}{k^{1+\alpha}},\qquad k\to\infty.
\end{align*}

Also the regular behavior of eigenvalues naturally arises when considering a Gaussian random field of the tensor product type, see Subsection~\ref{subsec:tensor-exact}. The first result on logarithmic asymptotics for such fields was obtained by E.~Cs\'aki~\cite{Csa1982} who investigated the Brownian sheet $\mathbb{W}^d$  on the unit cube ${\cal Q}_d$. His result, obtained by direct application of the Sytaya theorem (Theorem~\ref{thm:sytaya}), reads as follows:\footnote{Also Cs\'aki considered the pinned Brownian sheet 
$$
\stackrel{\circ}{\mathbb{W}}{}^{\!d}(x_1,\dots,x_d)=\mathbb{W}^d(x_1,\dots,x_d)-x_1\dots x_d\mathbb{W}^d(1,\dots,1)
$$
and obtained the following relation
$$
\log \mathbb{P}\{\|\stackrel{\circ}{\mathbb{W}}{}^{\!d}\|_{2,{\cal Q}_d}\leq \ep\}\sim\log \mathbb{P}\{\|\mathbb{W}^d\|_{2,{\cal Q}_d}\leq \ep\}, \qquad \ep\to0,
$$ 
that is in concordance with the fact that the pinned Brownian sheet is a one-dimensional perturbation of $\mathbb{W}^d$.

 Notice that in~\cite[Theorem~3.9]{Fat2003} it is erroneously stated that Cs\'aki considered only the case $d=2$.
}
\begin{align*}
\log \mathbb{P}\{\|\mathbb{W}^d\|_{2,{\cal Q}_d}\leq \ep\}\sim -\frac {2^{2d-5}}{\pi^{2d-2}((d-1)!)^2}\,\ep^{-2}\log^{2d-2}(\ep^{-1})
 \,, \qquad \ep\to0.
\end{align*} 

The result of~\cite{Csa1982} was generalized by W.V.~Li~\cite[Example 2]{Li1992} on the base of his comparison principle. Namely, he considered the double sums~\eqref{eq:double-sum} 
for $p_1,p_2>1$ and $\delta_1,\delta_2>-1$ and showed that the logarithmic $L_2$-small ball asymptotics for ${\cal S}^{(2)}$ in the cases $p_1=p_2$ and $p_1\ne p_2$ is essentially different. Namely, we have as $\ep\to0$
\begin{align}
\label{eq:log-P-Brownian-sheet-Li}
\log\mathbb{P}\{{\cal S}^{(2)}\leq \ep^2\}
\sim 
\begin{cases}
-\,C_1\,\ep^{-\frac{2}{p_1-1}}\log^{\frac{p_1}{p_1-1}}(\ep^{-1}), & p_1=p_2;
\\
-\,C_2\,\ep^{-\frac{2}{p_1-1}}, & p_1<p_2,
\end{cases}
\end{align}
where $C_1$ and $C_2$ are explicit constants; moreover, $C_1$ depends on $p_1$ only, whereas $C_2$ depends on $p_1$, $p_2$ and $\delta_2$ (but not on $\delta_1$!).

In~\cite[Example 3]{Li1992}, a generalization of~\eqref{eq:log-P-Brownian-sheet-Li} for multiple sums was given in the case of equal marginal exponents $p_i$.
\medskip

These results show that even in the case of purely power eigenvalues asymptotics of marginal processes, the asymptotics for the tensor product contains a slowly varying (logarithmic) factor. 

In the paper~\cite{KarNazNik2008-Tensor} A.I.~Karol', A.I.~Nazarov, and Ya.Yu.~Nikitin generalized the problem in a natural way and considered the fields-products with regularly varying marginal eigenvalues asymptotics of general form. They derived one-term asymptotics as $\tau\to0$ of the counting function ${\cal N}(\tau)$ for the eigenvalues of the tensor product 
in terms of spectra of marginal operators. This provides the logarithmic $L_2$-small ball asymptotics for corresponding field-product via Lemma~\ref{lm:log-P-reg-spect-asymp}.

It is sufficient to consider the case of two marginal operators. Let ${\cal K}_1$ and ${\cal K}_2$ be positive (compact, self-adjoint) operators. Suppose that
\begin{align}
\label{eq:NK1-NK2}
    {\cal N}_{{\cal K}_1}(\tau)\sim \tau^{-\frac 1{p_1}}\Phi_1(\tau^{-1}), \qquad {\cal N}_{{\cal K}_2}(\tau)\sim \tau^{-\frac1{p_2}}\Phi_2(\tau^{-1}), \quad \tau\to0,
\end{align}
(recall that ${\cal N}_{\cal K}(\tau)$ stands for the counting function of the eigenvalues of~${\cal K}$). Here the marginal exponents satisfy the assumption $p_1,p_2>1$, and the functions $\Phi_1$ and $\Phi_2$ are slowly varying at infinity. Also denote by
$$
{\cal N}_{\otimes}(\tau):={\cal N}_{{\cal K}_1\otimes{\cal K}_2}(\tau)
$$
the eigenvalues counting function of the tensor product.

It turns out that in general situation, as in~\eqref{eq:log-P-Brownian-sheet-Li}, two essentially different cases arise. 
If $p_1\ne p_2$ then the order of growth of ${\cal N}_{\otimes}(\tau)$ as $\tau\to0$ is defined by the ``factor'' with the slowest marginal exponent. However, the logarithmic small ball constant in this case is expressed in terms of zeta-function of the ``factor'' with fast eigenvalues decrease and requires complete information about its spectrum. 

In contrast, for the equal marginal exponents $p_1=p_2=p$, the order of growth of ${\cal N}_{\otimes}(\tau)$ as $\tau\to0$ deeply depends on the behavior of $\Phi_1$ and $\Phi_2$ at infinity. In this relation, the so-called \textit{\textbf{Mellin convolution}} of two SVFs 
\begin{align}
\label{eq:Mellin-conv}
(\Phi*\Psi)(s)
=
\int\limits_1^{s}\Phi(\sigma)\Psi
(s\sigma^{-1})\sigma^{-1}d\sigma
\end{align}
was investigated in~\cite[Section 2]{KarNazNik2008-Tensor}. In particular, it was shown that $\Phi*\Psi$ is also a SVF.

\begin{theorem}[{\cite[Theorems~3.2 and 3.4]{KarNazNik2008-Tensor}}]
\label{thm:count-function-tensor-prod}
Suppose that the operator ${\cal K}_1$ satisfies~\eqref{eq:NK1-NK2}.
\smallskip

{\bf 1}. Let ${\cal K}_2$  be an operator subject to
$$
{\cal N}_{{\cal K}_2}(\tau)= O(\tau^{-\frac 1{p_2}}),\quad \tau\to0, 
\qquad p_2>p_1.
$$
Then we have
$$
{\cal N}_{\otimes}(\tau) \sim C\,\tau^{-\frac 1{p_1}}\Phi_1(\tau^{-1}),\quad \tau\to0,
$$
where $C=\sum\limits_{k=1}^{\infty}
\mu_k^{\frac 1{p_1}}({\cal K}_2)$.
\smallskip

{\bf 2}. Suppose that the operator ${\cal K}_2$ satisfies~\eqref{eq:NK1-NK2}, and let $p_1=p_2=p$. Then we have
\begin{align}
\label{eq:N-tensor-prod-SVF}
{\cal N}_{\otimes}(\tau) \sim \tau^{-\frac 1p}{\Phi(\tau^{-1})},\qquad
\tau\to0, 
\end{align}
where $\Phi$ is explicit (though complicated) SVF depending on convergence of integrals\footnote{In particular, if both integrals~\eqref{eq:phi1-phi2} diverge then $\Phi(s)=\frac 1p(\Phi_1*\Phi_2)(s)$.}
\begin{align}
\label{eq:phi1-phi2}
\int\limits_1^{\infty}\Phi_1(\sigma)\sigma^{-1}d\sigma\qquad\mbox{and}\qquad \int\limits_1^{\infty}\Phi_2(\sigma)\sigma^{-1}d\sigma.
\end{align} 
\end{theorem}

As an illustration to these general results, we give the exact expression for the SVF $\Phi$ from~\eqref{eq:N-tensor-prod-SVF} in the case where SVFs in marginal asymptotics are powers of logarithm, see~\cite[Example 3]{KarNazNik2008-Tensor}. Namely, let $p_1=p_2=p$, and let
$$
\Phi_1(s)=\log^{\varkappa_1}(s),\qquad \Phi_2(s)=\log^{\varkappa_2}(s).
$$ 
Then the asymptotic formula~\eqref{eq:N-tensor-prod-SVF} holds with
$$
\Phi(s)=
\begin{cases}
\frac 1{p_{\displaystyle\vphantom 1}}\cdot{\bf B}(\varkappa_1+1,\varkappa_2+1)\log
^{\varkappa_1+\varkappa_2+1}(s),& \varkappa_1>-1,\varkappa_2>-1;
\\
\frac 1{p_{\displaystyle\vphantom 1}}\cdot\log(\log(s))\log^{\varkappa_2}(s),&
\varkappa_1=-1,\varkappa_2>-1;
\\ 
\frac 2{p_{\displaystyle\vphantom 1}}\cdot
\log(\log(s))\log^{-1}(s),& \varkappa_1=\varkappa_2=-1;
\\ 
 \sum\limits_{k=1_{\displaystyle\vphantom 1}}^{\infty}\mu_k^{\frac 1p}({\cal K}_1)\cdot\log^{\varkappa_2}(s),&
\varkappa_1<-1,\varkappa_1<\varkappa_2;
\\
\Bigl(\,\sum\limits_{k=1}^{\infty}\mu_k^{\frac 1p}({\cal K}_1) + \sum\limits_{k=1}^{\infty}\mu_k^{\frac 1p}({\cal K}_2)\Bigr)\cdot\log^{\varkappa_1}(s),&
\varkappa_1=\varkappa_2<-1
\end{cases}
$$
(here $\bf B$ is the Euler Beta-function).

The particular case $\varkappa_1=\varkappa_2=0$ of this formula can be extracted from~\cite{PapWas1990}, see also~\cite{LusPag2004}. For $\varkappa_1,\varkappa_2\ge0$, this result was obtained in~\cite[Section~3]{GaoLi2007-log}, see also~\cite{GaoLi2007-slepian}.

Using general results of Theorem~\ref{thm:count-function-tensor-prod} and Lemma~\ref{lm:log-P-reg-spect-asymp}, the authors of~\cite{KarNazNik2008-Tensor} considered many examples of fields-products (with equal or distinct marginal exponents) and sums of such products. Moreover, marginal processes can in turn be multiparameter, and a wide class of weighted $L_2$-norms is allowed.\footnote{More examples can be found in~\cite{Pus2008-Matern}, \cite{Pus2011-Matern}.}\medskip

More examples of regular eigenvalues asymptotics are generated by fractional Gaussian processes with variable Hurst index. The first example of such processes was apparently the \textit{\textbf{multifractional Brownian motion}}. It was introduced in \cite{PelVeh1995}, \cite{BenCohIst1998} and was investigated in various directions in several papers. There are some different definitions of this process equivalent up to a multiplicative deterministic function. We give the so-called harmonizable representation \cite{BenCohIst1998}
\begin{equation}\label{eq:mBm}
W^{H(\cdot)}(t)=C_*(H(t)) \int\limits_{-\infty}^{\infty}\frac{\exp(it\xi)-1}{|\xi|^{H(t)+\frac 12}}\,dW(\xi),
\end{equation} 
where $W(\xi)$ is a conventional Wiener process and the functional Hurst parameter $H(t)$ satisfies $0<H(t)<1$. The choice of normalizing factor
\begin{equation*}
C_*(H)= \left( \frac{\Gamma(2H +1)\sin(\pi H)}{2\pi}\right) ^{\frac 12}
\end{equation*}
ensures that the variance of $W^{H(\cdot)}(1)$ equals one.

A different process of the similar structure is the so-called \textit{\textbf{multifractal Brownian motion}} introduced in \cite{RalShe2009}, see also \cite{Ryv2015}.\footnote{For this process the functional Hurst parameter $H(t)$ satisfies $\frac 12<H(t)<1$.} 

Both mentioned processes have zero mean. For $H(x)\equiv H=const$ they both coincide with conventional fBm.
\medskip

The $L_2$-small ball asymptotics for these processes was derived by A.I.~Ka\-rol' and A.I.~Nazarov~\cite{KarNaz2021-Multifrac}, \cite{Kar2022} under some mild assumptions on the regularity of $H(t)$. To obtain corresponding eigenvalues asymptotics, authors used the approximate spectral projector method in combination with the methods of asymptotic perturbation theory.

It turns out that the asymptotic properties of the eigenvalues are determined by the set where the function $H(t)$ attains its minimal value $H_{\min}$. If the measure of this set is positive, the counting function ${\cal N}(\tau)$ has classical power asymptotics. The situation where this measure vanishes is much more complicated. In this case, the asymptotics is non-power (regularly varying at zero) and is determined by the behavior as $\epsilon\to0$ of the function 
$$
{\cal V}(\epsilon):=\big|\{t\in[0,1] : H_{\min}< H(t)<H_{\min}+\epsilon\}\big|.
$$

We restrict ourselves to three examples from~\cite{KarNaz2021-Multifrac}, \cite{Kar2022}. For simplicity we assume that $H_{min}=\frac 12$. 
\medskip

{\bf 1}. Let $H(t)=\frac 12+(t-t_0)_+^\gamma$, $0<t_0\le 1$, $\gamma>0$. Then the minimal value set is $[0,t_0]$, and Theorem~2 in~\cite{KarNaz2021-Multifrac} gives
$$
\log \mathbb{P}\{\|W^{H(\cdot)}\|_{2,[0,1]}\leq \ep\}\sim-\,\frac{t_0^2}8\, \ep^{-2}, \qquad \ep\to0.
$$
For $t_0=1$ we deal with standard Wiener process on $[0,1]$, and this result is well known. 
\medskip

{\bf 2}. Let $H(t) = \frac 12 + |t-t_0|^\gamma$, $\gamma>0$. In this case we have
$$
{\cal V}(\epsilon)= 
\begin{cases}
    2\epsilon^{\frac 1{\gamma}}, & 0<t_0<1;
    \\
    \epsilon^{\frac 1{\gamma}}, & t_0=0,1.
\end{cases}
$$
Therefore, Theorem~2 in~\cite{KarNaz2021-Multifrac} gives
\begin{equation*}
\log \mathbb{P}\{\|W^{H(\cdot)}\|_{2,[0,1]}\leq \ep\}\sim-\, C(t_0)\Gamma^2\big(1+\tfrac 1\gamma\big) \ep^{-2} \log^{-\frac 2\gamma} (\ep^{-1}), \quad \ep\to0,
\end{equation*}
where
$C(t_0)=2^{-1-\frac 2{\gamma}}$ if $0<t_0<1$ and $C(t_0)=2^{-3-\frac 2{\gamma}}$ if $t_0=0,1$.
\medskip

{\bf 3}. Let $H(t)=\frac 12 +{\rm dist}^\gamma(t, \mathfrak{D})$, where $\mathfrak{D}$ is standard Cantor set. A tedious but simple calculation gives for small $\epsilon>0$
$$
{\cal V}(\epsilon)= \epsilon^{1-\frac{\log (2)}{\log (3)}}\,\zeta(\log (\epsilon^{-1})),
$$
where $\zeta$ is a  {\it periodic} function. In this case 
Theorem~2 in~\cite{KarNaz2021-Multifrac} is not applicable. However, more refined Theorem~4.4 in~\cite{Kar2022} gives
\begin{align}
\label{eq:log-P-mBm-Cantor-set}
\log \mathbb{P}\{\|W^{H(\cdot)}\|_{2,[0,1]}\leq \ep\}\sim \ep^{-2}\log^{-\frac 2{\gamma}}(\ep^{-1})\eta(\log(\log(\ep^{-1}))), \quad \ep\to0,
\end{align}
where $\eta$ is an explicit (though complicated) continuous positive $\gamma\log(3)$-pe\-ri\-od\-ic function.\footnote{We underline that the asymptotics in~\eqref{eq:log-P-mBm-Cantor-set} is still regularly varying.}
\medskip

Also it is shown in~\cite{KarNaz2021-Multifrac}, \cite{Kar2022} that, despite the fact that the behavior of covariances of multifractional and multifractal Brownian motions is significantly different, for $H(t) > \frac 12$ their logarithmic $L_2$-small ball asymptotics coincide.

\subsection{Almost regular asymptotics}
\label{subsec:alm-reg-asymp}
Here we consider processes with more complicated (the so-called \textit{\textbf{almost regular}}) eigenvalues asymptotics\footnote{Similarly to Footnote~\ref{foot:count-eigen}, the relation~\eqref{eq:almost-reg-eigen} is equivalent to the following asymptotics of the counting function:
$$
{\cal N}(\tau)\sim \tau^{-\frac 1p}\Phi(\tau^{-1})\mathfrak{g}(\log(\tau^{-1})),\quad \tau\to0,
$$
where $\Phi$ is a SVF and $\mathfrak{g}$ satisfies the Condition {\bf A}. Moreover, $\mathfrak{g}$ can be chosen $T$-periodic if $\mathfrak{h}$ is $\frac Tp$-periodic.
}
\begin{align}
\label{eq:almost-reg-eigen}    
\mu_k\sim \,k^{-p}\Psi(k)\mathfrak{h}(\log(k)), \quad k\to\infty,
\end{align}
where $p>1$, $\Psi$ is a SVF, and the function $\mathfrak{h}$ satisfies the following condition:\medskip

{\bf Condition A}. \textit{The function is positive, uniformly continuous on $\mathbb{R}_+$, bounded and separated from zero.}

\begin{remark}
    In general case the above mentioned conditions do not imply the monotonicity of the function in the right hand side of~\eqref{eq:almost-reg-eigen}. In what follows we always assume such monotonicity.
\end{remark}

We begin with the generalization of the relation~\eqref{eq:power-asymp-small-balls} for the almost regular asymptotics:

\begin{lemma}[{\cite[Theorem~8]{Rast2018-tensor}}]\!\!\footnote{ For $\Psi\equiv1$ this fact was obtained earlier in~\cite[Theorem~4.2]{Naz2004-logselfsim}.}
\label{lm:log-P-alm-reg-spect-asymp}
Let $p>1$. Suppose that the function $\Psi(t)$ is slowly varying at infinity, and the function $\mathfrak{h}$ satisfies the Condition {\bf A}. Then, as $\ep\to0$,
\begin{align}
\label{eq:alm-reg-asymp-small-balls}
\log\mathbb{P}\Big\{\sum_{k=1}^{\infty}k^{-p}\Psi(k)\mathfrak{h}(\log(k)) \,\xi_k^2 \leq
\ep^2\Big\}
\sim - \,\ep^{-\frac {2}{p-1}}\Theta(\ep^{-1})\zeta(\log (\ep^{-1})),
\end{align}
here $\Theta$ is the same SVF as in~\eqref{eq:reg-asymp-small-balls}, and the function $\zeta$ satisfies Condition~{\bf A}. Moreover, if $\mathfrak{h}$ is $\frac Tp$-periodic then $\zeta$ can be chosen $\frac {T(p-1)}{2p}$-periodic.

\end{lemma}

The spectral asymptotics of the form~\eqref{eq:almost-reg-eigen} (with $\Psi\equiv 1$) can arise if we deal with a Green Gaussian process $X$ in the space $L_2[0,1;d\mathfrak{m}]$, where $\mathfrak{m}$ is a measure of special class (see below). In this case $\mu_k$ are the eigenvalues of the following integral equation:
\begin{align}
    \label{eq:int-eq-eigen-meas}
    \int\limits_0^1 G_X(t,s)\varphi(s)\,d\mathfrak{m}(s)=\mu\varphi(t), \qquad 0\le t\le 1,
\end{align}
It was shown in~\cite{BirSol1970} (even in multidimensional case) that if the measure $\mathfrak{m}$ contains absolutely continuous component then its singular component does not influence upon the main term of spectral asymptotics. Therefore, the singular component of $\mathfrak{m}$ does not influence upon the logarithmic small ball asymptotics in $L_2[0,1;d\mathfrak{m}]$ for corresponding processes, see Theorem~\ref{thm:log-asymp-singularity-power}.

The situation drastically changes when $\mathfrak{m}$ is singular with respect to Lebesgue measure. V.V.~Borzov~\cite{Borz1971} showed that if $G_X$ is the Green 
function of the problem~\eqref{eq:L-bvp} for the operator ${\cal L}u\equiv(-1)^{\ell}u^{(2\ell)}$ then the corresponding
eigenvalues $\mu_k$ obey the estimate $o(k^{-2\ell})$ instead of usual asymptotics $\mu_k\sim C\cdot k^{-2\ell}$ in the 
case of non-singular $\mathfrak{m}$. For some special classes of $\mathfrak{m}$ better upper estimates were obtained in~\cite{Borz1971}. 

Further improvements of these results are possible if we restrict ourselves to the so-called \textit{\textbf{self-similar measures}} $\mathfrak{m}$.
\medskip

Recall the construction of self-similar probability measure on $[0,1]$ (general situation is described in~\cite{Hut1981-frac}). Consider $n\ge2$ 
non-empty non-in\-ter\-sect\-ing intervals $I_j\subset (0,1)$, $j=1,\dots,n$, and define a family of affine functions $S_j$, $j=1,\dots,n$, mapping $(0,1)$ onto $I_j$. Consider also positive numbers (weights) $\mathfrak{r}_j$, $j=1,\dots,n$, such that $\sum\limits_j\mathfrak{r}_j=1$.

It is not difficult to show (see, e.g., \cite[Section~2]{Naz2004-logselfsim}), that there is a unique probability measure $\mathfrak{m}$ such that for any 
Lebesgue-measurable set ${\cal E}\subset[0,1]$ 
\begin{align}
\label{eq:self-similar-meas}
\mathfrak{m}({\cal E})= \sum\limits_{j=1}^n\mathfrak{r}_j\cdot\mathfrak{m}(S_j^{-1}({\cal E})).
\end{align}
The relation~\eqref{eq:self-similar-meas} shows the \textit{\textbf{self-similarity}} property of the measure $\mathfrak{m}$. Notice that $\mathfrak{m}$ has no atoms. Its primitive is called a \textit{\textbf{generalized Cantor ladder}}.

When $\sum\limits_j|I_j|<1$, the support of $\mathfrak{m}$ (minimal {\it closed} set 
${\cal E}\subset[0,1]$ such that $\mathfrak{m}([0,1]\backslash{\cal E})=0$) is called 
a \textit{\textbf{generalized Cantor set}}. Its Hausdorff 
dimension $\alpha\in(0,1)$ is equal to the unique solution of the equation
$$
\sum\limits_{j=1}^n|I_j|^{\alpha}=1.
$$
In the case $\sum\limits_j|I_j|=1$ the support of $\mathfrak{m}$ is $[0,1]$, 
and $\alpha=1$. If, in addition, $\mathfrak{r}_j=|I_j|$, $j=1,\dots,n$, then $\mathfrak{m}$ is the conventional Lebesgue measure. However in other cases $\mathfrak{m}$ is singular.

\begin{figure}
    \centering
    \includegraphics[width=0.4\textwidth]{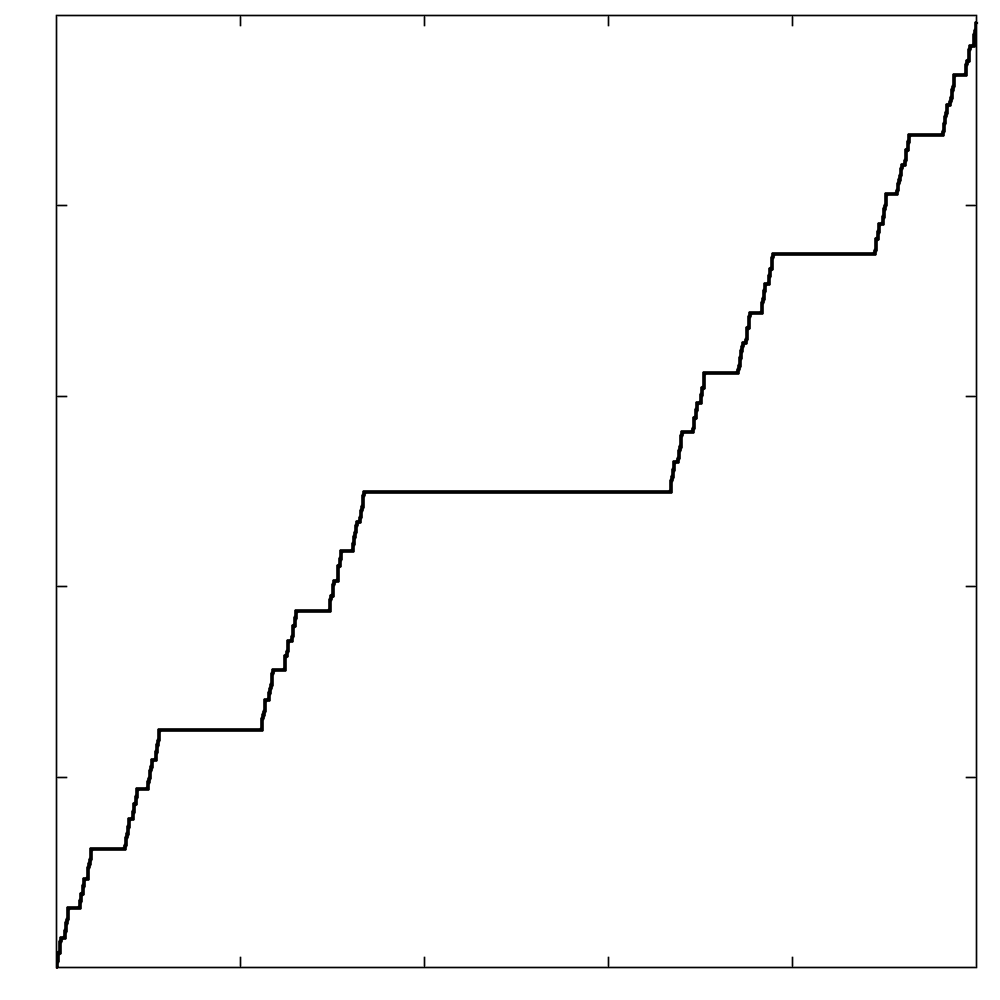}
    \caption{Cantor ladder}
    \label{fig:cantor-ladder}
\end{figure}

In the case $n=2$, $I_1=(0,\frac 13)$, $I_2=(\frac 23,1)$, $\mathfrak{r}_1=\mathfrak{r}_2=\frac 12$, $\mathfrak{m}$ is the standard Cantor measure, and its primitive is the standard Cantor ladder, see Fig.~\ref{fig:cantor-ladder}.

Recall that we deal with a Green
Gaussian process, so the integral equation~\eqref{eq:int-eq-eigen-meas} is equivalent to the boundary value problem\footnote{The equation in~\eqref{eq:L-meas} is understood in the sense of distributions.}
\begin{align}
\label{eq:L-meas}
{\cal L}u=\lambda \mathfrak{m}u \quad {\rm on} 
\quad [0,1], \qquad u\in {\cal D}({\cal L}),
\end{align}
where $\lambda=\mu^{-1}$, ${\cal L}$ is defined in~\eqref{eq:L}, and ${\cal D}({\cal L})$ is determined by the boundary conditions~\eqref{eq:bc-non-separated}. For simplicity we assume that $p_\ell\equiv1$.

We list the main steps in the research on the problem~\eqref{eq:L-meas}. More references can be found in~\cite{Naz2004-logselfsim}. 

For the simplest Sturm--Liouville operator ${\cal L}u\equiv-u''$, the exact power order $p$ of the eigenvalues growth was obtained by T.~Fujita~\cite{Fuj1987}. The one-term eigenvalues asymptotics in this case was obtained independently by J.~Kigami and M.L.~Lapidus~\cite{KigLap1993} and by M.Z.~Solomyak and E.~Verbitsky~\cite{SolVer1995}. A.I.~Nazarov~\cite{Naz2004-logselfsim} generalized their result for the operators of arbitrary order.

\begin{theorem}[{\cite[Theorem~3.1]{Naz2004-logselfsim}}]
\label{thm:self-sim-asymp-count-funct}
Given self-similar probability measure~$\mathfrak{m}$, define
$$
\mathfrak{c}_j:=\mathfrak{r}_j\cdot|I_j|^{2\ell-1}, \quad j=1,\dots,n,
$$

\noindent and denote by $p$ the unique solution of the equation
$$
\sum\limits_{j=1}^n \mathfrak{c}_j^{\frac1p}=1.
$$
Evidently, $p\ge2\ell$. Moreover, $p=2\ell$ only if $\mathfrak{m}$ is the Lebesgue measure.

Let $X$ be a Green Gaussian process, corresponding to the order~$2\ell$ ODO.
\medskip

{\bf 1}. If at least one 
ratio $\log(\mathfrak{c}_i^{-1})/\log(\mathfrak{c}_j^{-1})$ is irrational, 
the self-similarity of the measure $\mathfrak{m}$ is called \textit{\textbf{non-arithmetic}}. In this case the counting function of eigenvalues to the integral equation~\eqref{eq:int-eq-eigen-meas} has purely power asymptotics:
$$
{\cal N}(\tau)\sim M\tau^{-\frac 1p}, 
\quad \tau\to0,
$$
where $M$ is some (unknown) constant.
\medskip

{\bf 2}. If all quantities $\log(\mathfrak{c}_j^{-1})$ are 
mutually commensurable, the self-sim\-i\-lar\-i\-ty  of $\mathfrak{m}$ is called \textit{\textbf{arithmetic}}. In this case 
\begin{align}
    \label{eq:N-asymp-self-sim-arithm}
    {\cal N}(\tau)\sim \tau^{-\frac1p}\mathfrak{g}(\log(\tau^{-1})), 
\quad \tau\to0,
\end{align} 
where the function $\mathfrak{g}$ satisfies Condition~{\bf A}. Moreover, $\mathfrak{g}$ is periodic, and its period is equal to the greatest common 
divisor of $\log(\mathfrak{c}_j^{-1})$, $j=1,\dots,n$.
\end{theorem}

\begin{remark}
\label{rm:spect-dim}
 The quantity\footnote{Lemma~\ref{lm:log-P-alm-reg-spect-asymp} implies that the power order of the logarithmic $L_2$-small ball asymptotics for the corresponding Green Gaussian process is just $\frac {2\gamma}{2\ell-1}$.}
$\gamma=\frac {2\ell-1}{p-1}$
is called \textit{\textbf{spectral dimension of order 
$2\ell-1$}} of the self-similar measure~$\mathfrak{m}$. It is shown in~\cite[Theorem~5.3]{Naz2004-logselfsim} that $\gamma\le\alpha$ (recall that $\alpha$ is the Hausdorff dimension of the support of $\mathfrak{m}$), and $\gamma=\alpha$ if and only if $\mathfrak{r}_j=|I_j|^{\alpha}$ 
for any $j=1,\dots,n$. In particular, this is the case for the standard Cantor measure, where $\gamma=\alpha=\frac {\log(2)}{\log(3)}$.
\end{remark}

\begin{remark}
Theorem~\ref{thm:self-sim-asymp-count-funct} was generalized in several directions. Namely, the problem~\eqref{eq:L-meas} was considered under the following assumptions:
\begin{enumerate}
\item ${\cal L}u\equiv -u''$, and $\mathfrak{m}$ is the distributional derivative of the Minkowski question-mark function,\footnote{The function $?(t)$ introduced by H.~Minkowski~\cite{Min1904} is well known in number theory.} see~\cite{She2015}. In this case $\mathfrak{m}$ is a non-affine self-similar measure, i.e. relation~\eqref{eq:self-similar-meas} holds with non-affine diffeomorphisms $S_j$, namely,
$$
n=2, \quad S_1(t)=\frac t{t+1},\quad S_2(t)=\frac 1{2-t}. 
$$
In~\cite[Theorem~3.2]{She2015}, two-sided estimates for the eigenvalues of counting function were obtained:
$$
{\cal N}(\tau)\asymp \tau^{\frac {\mathfrak{d}}{\mathfrak{d}+1}},\quad \tau\to0; \qquad \mathfrak{d}=\frac 12 \Big(\int\limits_0^1\log_2(1+t)\,d\mathfrak{m}\Big)^{-1}.
$$

\item ${\cal L}u\equiv -u''$, and $\mathfrak{m}$ is a \textit{\textbf{self-conformal}} measure, i.e. $S_j$ in the relation~\eqref{eq:self-similar-meas} are non-affine diffeomorphisms mapping $(0,1)$ onto $I_j$, subject to some additional assumptions, see~\cite{FreRast2020}, \cite{FreRast2020-sib};
\item $\mathfrak{m}$ is a sign-changing self-similar distribution,  see~\cite{VlaShe2006-indef}, \cite{VlaShe2006};
\item ${\cal L}$ is a so-called \textit{\textbf{measure geometric Krein--Feller operator}} related to two self-similar measures, see~\cite{Fre2011} and references therein.\footnote{Notice, however, that some of proofs in~\cite{Fre2011} contain gaps, see~\cite[2.2.3]{Vla2019}.}
\end{enumerate}
Up to our knowledge, the generalizations mentioned in items 3 and 4 have not found an application in the study of small ball probabilities yet.
\end{remark}

Using general results of Theorem~\ref{thm:self-sim-asymp-count-funct} and Lemma~\ref{lm:log-P-alm-reg-spect-asymp}, the logarithmic small ball asymptotics in $L_2[0,1;d\mathfrak{m}]$ were obtained in~\cite{Naz2004-logselfsim} for several classical Green Gaussian processes.
\medskip

The statement of ``arithmetic'' part {\bf 2} in Theorem~\ref{thm:self-sim-asymp-count-funct} does not exclude the ``degenerate'' case where function $\mathfrak{g}$ in~\eqref{eq:N-asymp-self-sim-arithm} is a constant, i.e. 
${\cal N}(\tau)$ has classical power asymptotics. It was conjectured in~\cite{Naz2004-logselfsim} that it is not the case, i.e. $\mathfrak{g}\ne const$ for any non-Lebesgue arithmetically self-similar measure~$\mathfrak{m}$. 

A.A.~Vladimirov and I.A.~Sheipak~\cite{VlaShe2006} confirmed this conjecture for the standard Cantor measure and ${\cal L}u\equiv-u''$ using computer-assisted proof. Later~\cite{VlaShe2013} they succeeded to give a complete characterization of the function $\mathfrak{g}$ under assumption that $\mathfrak{m}$ is a so-called \textit{\textbf{even}} Cantor-type measure.

\begin{theorem}[{\cite[Section~5]{VlaShe2013}}] 
\label{thm:almost-reg-g}
Suppose that $\mathfrak{m}$ is an even Cantor-type measure, that is, for all $j=1,\dots,n$, the lengths $|I_j|$ are equal, the gaps between neighbor intervals $I_j$ are equal (and non-zero), and the weights~$\mathfrak{r}_j$ are equal. Let ${\cal L}u\equiv-u''$. Then 
the function $\mathfrak{g}$ in~\eqref{eq:N-asymp-self-sim-arithm} satisfies the relation 
\begin{align}
\label{eq:g-non-const}
\mathfrak{g}(t)=\exp(-\tfrac tp)\mathfrak{s}(t),
\end{align}
where $\mathfrak{s}$ is some purely singular non-decreasing function (i.e. a primitive of a singular measure).\footnote{Formula~\eqref{eq:g-non-const} evidently implies that $\mathfrak{g}\ne const$.}
\end{theorem}

In~\cite{Vla2016} this result was transferred to the fourth order operator ${\cal L}$. 

N.V. Rastegaev~\cite{Rast2015-cantor}, \cite{Rast2020-selfsimilar} obtained formula~\eqref{eq:g-non-const} for ${\cal L}u\equiv-u''$ and arbitrary arithmetically self-similar measure $\mathfrak{m}$ under the only assumption: all gaps between neighbor intervals $I_j$ are non-zero. Further generalization (in particular, for the higher order operator ${\cal L}$) is an interesting open problem.
\medskip

All previous examples generate the eigenvalues asymptotics~\eqref{eq:almost-reg-eigen} with $\Psi\equiv1$. A general form of~\eqref{eq:almost-reg-eigen} arises, for instance,  when considering tensor products. We describe briefly corresponding results obtained by N.V. Rastegaev~\cite{Rast2018-tensor}.
Let ${\cal K}_1$ and ${\cal K}_2$ be positive (compact, self-adjoint) operators. Suppose that 
\begin{align}
\label{eq:NK1-almost}
    {\cal N}_{{\cal K}_1}(\tau)\sim \tau^{-\frac 1{p_1}}\Phi_1(\tau^{-1})\mathfrak{g}_1(\log(\tau^{-1})), \quad  
    \tau\to0,
\end{align}
where $p_1>1$, $\Phi_1$ is a SVF, and $\mathfrak{g}_1$ is a continuous positive $T_1$-periodic function.\footnote{The relation~\eqref{eq:NK1-almost} implies that the function $\mathfrak{g}$ has bounded variation on $[0,T_1]$, see~\cite[Lemma~6]{Rast2018-tensor}.}

If ${\cal N}_{{\cal K}_2}(\tau)= O(\tau^{-\frac 1{p_2}})$ as $\tau\to0$, with $p_2>p_1$ then
$$
{\cal N}_{\otimes}(\tau) \sim \tau^{-\frac 1{p_1}}\Phi_1(\tau^{-1})\mathfrak{g}^*(\log(\tau^{-1})),\quad \tau\to0,
$$
where $\mathfrak{g}^*$ is a continuous positive $T_1$-periodic function. 

\begin{remark}
It is shown in~\cite[Remark~2]{Rast2018-tensor} that in general case $\mathfrak{g}^*$ can degenerate into constant even if $\mathfrak{g}_1$ is not a constant. However, if $\mathfrak{g}_1$ has the form given in Theorem~\ref{thm:almost-reg-g} then $\mathfrak{g}^*$ also has such a form.
\end{remark}

Suppose now that ${\cal K}_2$ satisfies a similar relation 
$$
{\cal N}_{{\cal K}_2}(\tau)\sim \tau^{-\frac1{p_2}}\Phi_2(\tau^{-1})\mathfrak{g}_2(\log(\tau^{-1})), \quad  \tau\to0
$$
($p_2>1$, $\Phi_2$ is a SVF, and $\mathfrak{g}_2$ is a continuous positive $T_2$-periodic function),
and let $p_1=p_2=p$. Then the asymptotics of ${\cal N}_{\otimes}(\tau)$ as $\tau\to0$ is more complicated and depends on convergence of integrals~\eqref{eq:phi1-phi2} and, if they both diverge, on the commensurability of $T_1$ and $T_2$.
In particular, if both integrals~\eqref{eq:phi1-phi2} diverge and $T_1=T_2=T$ then
$$
{\cal N}_{\otimes}(\tau) \sim \tau^{-\frac 1p}(\Phi_1*\Phi_2)(\tau^{-1})\mathfrak{g}_{\otimes}(\log(\tau^{-1})),\quad \tau\to0,
$$
(the Mellin convolution of SVFs is defined in~\eqref{eq:Mellin-conv}), where $\mathfrak{g}_{\otimes}$ is an explicit (though complicated) continuous positive $T$-periodic function. 

\begin{remark}    
The question of non-constancy of $\mathfrak{g}_{\otimes}$ remains open even in the case where both $\mathfrak{g}_1$ and $\mathfrak{g}_2$ have the form given in Theorem~\ref{thm:almost-reg-g}.
\end{remark}

\subsection{Slowly varying asymptotics}

In this subsection we consider processes with super-power decay of eigenvalues. More exactly, we assume that the counting function ${\cal N}(\tau)$ is slowly varying at zero. Corresponding analog of the relation~\eqref{eq:power-asymp-small-balls} reads as follows:

\begin{lemma}[{\cite[Theorem~2]{Naz2009-log}}]
\label{lm:log-P-slow-var-spect-asymp}
Let $(\mu_k)$, $k\in\mathbb N$, be a
positive sequence with counting function ${\cal N}(\tau)$, and let
\begin{align}
\label{eq:N-SVF}
    {\cal N}(\tau)\sim \Phi(\tau^{-1}), \quad \tau\to 0,
\end{align}
where $\Phi$ is a SVF. Then, as $\ep\to0$,
$$
\log\mathbb{P}\Big\{\sum_{k=1}^\infty\mu_k\xi_k^2\le \ep^2
\Big\} \sim-\frac {\Theta(\gamma)}2:=-\frac {(\Phi*{\bf 1})(\gamma)}2\equiv -\frac 12 \int\limits
_1^\gamma \Phi(\sigma)\sigma^{-1}d\sigma,
$$

\noindent where $\gamma=\gamma(r)$ satisfies the relation
$$
\frac {\Phi(\gamma)}{2\gamma}\sim \ep^2, \qquad \ep\to 0.
$$
\end{lemma}

\begin{remark}
    As it was mentioned in Subsection~\ref{subsec:reg-var}, $\Theta$ is a SVF, cf.~\eqref{eq:Mellin-conv}.
\end{remark}

The first examples of Gaussian processes satisfying~\eqref{eq:N-SVF} arose in~\cite{GaoLi2007-log}, where the case $\log(\mu_k)\sim -k^{\alpha}\Psi(k)$ as $k\to\infty$ with a SVF $\Psi$ was considered.\footnote{We notice that the result of~\cite[Theorem~2.2]{GaoLi2007-log} 
$$
\log\mathbb{P}\Big\{\sum_{k=1}^\infty\mu_k\xi_k^2\le \ep^2
\Big\} \sim -\frac {\alpha 2^{\frac 1\alpha}}{\alpha+1} \frac {\log^{\frac 1\alpha+1}(\ep^{-1})}{\Psi^{\frac 1\alpha}(\log^{\frac 1\alpha}(\ep^{-1}))},\quad \ep\to0,
$$
is not true in general case. However, it holds if, say, $\Psi(k)=\log^{\sigma}(k+1)$, $\sigma\in\mathbb{R}$.
}

The relation~\eqref{eq:N-SVF} is typical for the processes with smooth covariances. 
In~\cite{Naz2009-log}, a set of smooth stationary Gaussian processes ${\cal R}_{C,\alpha}$ ($C,\alpha>0$) was considered. These processes have zero mean-value and the
spectral density
\begin{align}
\label{eq:spect-meas-smooth}
    {\bf m}_{C,\alpha}(\xi)=\exp(-C|\xi|^{\alpha}),
\qquad \xi\in\mathbb R.
\end{align}
For instance, it is well known that
$$G_{{\cal R}_{C,1}}(s,t)=\frac C{\pi(C^2+(s-t)^2)}, \qquad
G_{{\cal R}_{C,2}}(s,t)=\frac 1{2\sqrt{\pi C}}\exp\left(-\frac
{(s-t)^2}{4C}\right).
$$

The spectral asymptotics for the integral operators with kernels of this type was treated in remarkable paper
\cite{Wid1964}. The following relation for the counting function can be extracted from~\cite[Theorems 1 and 2]{Wid1964}: as $\tau\to0$, we have
$$
{\cal N}(\tau)\sim\Phi(\tau^{-1})\equiv
\begin{cases}
\frac 1{\pi_{\displaystyle\vphantom 1} C^{\frac 1\alpha}} \cdot\log^{\frac 1\alpha}(\tau^{-1}), &
\ \alpha<1;\\ 
\frac 1{\pi_{\displaystyle\vphantom 1} \mathfrak {C}} \cdot\log(\tau^{-1}), &
\ \alpha=1;\\ 
\frac {\alpha}{2\alpha-2} \cdot\frac {\log(\tau^{-1})}{\log(\log(\tau^{-1}))}, &
\ \alpha>1.
\end{cases}
$$
Here
$$
\mathfrak{C}=\frac {{\bf K}\bigl(\mbox{sech}\bigl(\frac{\pi}{2C}\bigr)\bigr)}{{\bf K}
\bigl(\tanh\bigl(\frac{\pi}{2C}\bigr)\bigr)},
$$
while $\bf K$ is the complete elliptic
integral of the first kind.\medskip

So, Lemma~\ref{lm:log-P-slow-var-spect-asymp} yields the following relation as $\ep\to0$:
$$
\log{\mathbb P}\left\{\|{\cal R}_{C,\alpha}\|_{2,[0,1]}\le\ep\right\} 
\sim
\begin{cases}
- \,\big(\frac 2C\big)^{\frac
1\alpha}\, \frac {\alpha}{(\alpha+1)\pi_{\displaystyle\vphantom 1}}\cdot \log^{\frac {\alpha+1}\alpha} (\ep^{-1}),&\ \alpha<1;
\\
- \,
\frac 1{\pi_{\displaystyle\vphantom 1}\mathfrak{C}}\cdot\log^2(\ep^{-1}), &\ \alpha=1;
\\
- \,
\frac {\alpha}{2\alpha- 2} \cdot\frac {\log^2(\ep^{-1})}
{\log(\log(\ep^{-1}))},&\ \alpha>1.
\end{cases}
$$
\begin{remark}
Notice that the order of logarithmic $L_2$-small ball asymptotics in the latter case depends neither on $C$ nor even on $\alpha$. The same effect was found in~\cite{AurIbrLifZan2009} for the sup-norm.
\end{remark}

F.~Aurzada, F.~Gao, Th.~K{\"u}hn, W.V.~Li and Q.-M.~Shao  \cite{AurGaoKuhnLiShao2013} investigated the small ball behavior in $L_2$- and sup-norms of smooth self-similar Gaussian processes $X_{\alpha,\beta}$ with the covariance function
\begin{equation}
\label{eq:Laptev}
G_{X_{\alpha,\beta}}(t,s)=\frac {t^\alpha s^\alpha}{(t+s)^{2\beta+1}}.
\end{equation}
In particular, it was shown that under natural assumption $\alpha>\beta>-\frac 12$ we have
\begin{align}
\label{eq:log-P-smooth-5authors}
    \log{\mathbb P}\left\{\|X_{\alpha,\beta}\|_{2,[0,1]}\le\ep\right\} 
\sim -\,\frac 1{3(\alpha-\beta)\pi^2}\,\log^3(\ep^{-1}),\quad \ep\to0.
\end{align}

\begin{remark}
The method of \cite{AurGaoKuhnLiShao2013} is based on the estimates of entropy numbers of a linear operator generating the process, see~\cite{KueLi1993}.
We notice that the spectral asymptotics for the operator with kernel~\eqref{eq:Laptev} was obtained many years ago by A.A.~Laptev~\cite{Lap1974}. So, formula~\eqref{eq:log-P-smooth-5authors} can be obtained from the result of~\cite{Lap1974} and Lemma~\ref{lm:log-P-slow-var-spect-asymp}.
\end{remark}

Another class of problems where the relation~\eqref{eq:N-SVF} appears is related to the Green Gaussian process in the space $L_2[0,1;d\mathfrak{m}]$ with very ``poor'' measure~$\mathfrak{m}$.
We describe here (not in a full generality) corresponding class of \textit{\textbf{degenerate self-similar measures}}\footnote{This measure is not self-similar in the sense of~\cite{Hut1981-frac}.} on $[0,1]$, see~\cite{She2010}, \cite{NazShei2012}.

Let $0<t_1<\ldots <t_n<1$, $n\ge2$, be a partition of the segment $[0,1]$. We select one of intervals between these points $I=(t_k,t_{k+1})$ and define the affine function $S$ mapping $(0,1)$ onto $I$. Consider also positive numbers (weights) $a_j>0$, $j=1,\ldots,n$, and $\mathfrak{r}>0$ such that $\sum\limits_j a_j+\mathfrak{r}=1$.
 
It is not difficult to show (see, e.g., \cite{She2010}), that there is a unique probability measure $\mathfrak{m}$ such that for any Lebesgue-measurable set ${\cal E}\subset[0,1]$ 
\begin{align}
\label{eq:self-similar-degen-meas}
\mathfrak{m}({\cal E})= \mathfrak{r}\cdot\mathfrak{m}(S^{-1}({\cal E})) + \sum\limits_{t_j\in{\cal E}} a_j.
\end{align}

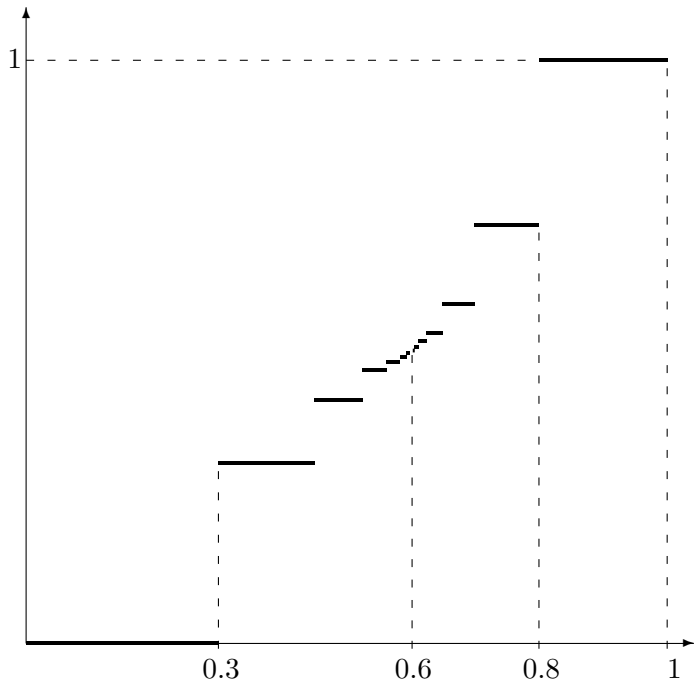
\begin{figure}[ht]
\begin{center}
 \hfil
 \begin{picture}(300,250)    
\put(0,20){\vector(1,0){250}}
\put(0,20){\vector(0,1){240}}

\linethickness{0.5mm}
\put(0,20){\line(1,0){72}}
\put(240,240){\line(-1,0){48}}

\put(72,88){\line(1,0){36}}
\put(192,178){\line(-1,0){24}}

\put(108,112){\line(1,0){18}}
\put(168,148){\line(-1,0){12}}

\put(126,123){\line(1,0){9}}
\put(156,137){\line(-1,0){6}}

\put(135,126){\line(1,0){5}}
\put(150,134){\line(-1,0){3}}

\put(140,128){\line(1,0){2.5}}
\put(147,132){\line(-1,0){1.5}}

\put(142.5,129.5){\line(1,0){1.25}}
\put(145.5,130.5){\line(-1,0){0.75}}

\thinlines
\multiput(0,240)(8,0){24}{\line(1,0){3}}
\put(-7,237){\small 1}
\put(240,7){\small 1}
\multiput(240,18)(0,8){28}{\line(0,1){3}}
\multiput(144.5,18)(0,7.7){15}{\line(0,1){3}}
\put(138,7){\small 0.6}

\multiput(72,18)(0,8){9}{\line(0,1){3}}\put(66,7){\small 0.3}
\multiput(192,18)(0,7.7){21}{\line(0,1){3}}\put(186,7){\small 0.8}
\end{picture}
\caption{The primitive of a degenerate self-similar measure} 
\label{fig:self-simil-degen-meas}
\end{center}
\end{figure}

In contrast with self-similar measures considered in Subsection~\ref{subsec:alm-reg-asymp}, the measure $\mathfrak{m}$, defined by relation~\eqref{eq:self-similar-degen-meas}, is discrete. Its support has a unique accumulation point $\widehat x$. Straightforward calculation shows that
\begin{equation}
\label{eq:osob_tochka}
\widehat x=\dfrac{t_k}{1-|I|}.
\end{equation}
The primitive of $\mathfrak{m}$ is a piecewise constant function. Fig.~\ref{fig:self-simil-degen-meas} shows its graph for the following values of parameters: $n=2$; $t_1=0.3$, $t_2=0.8$; $a_1=\frac 13$, $a_2=\frac 13$, $\mathfrak{r}=\frac 13$. 
Formula~(\ref{eq:osob_tochka}) gives $\hat x=0.6$.
\bigskip

\begin{remark}
Evidently, the Hausdorff dimension of the support of degenerate self-similar measure is equal to zero. Therefore, 
the spectral dimension of $\mathfrak{m}$ (see Remark~\ref{rm:spect-dim}) is also equal to zero. 
\end{remark}

The first result on the eigenvalues asymptotics for the boundary value problem~\eqref{eq:L-meas} with degenerate self-similar measure $\mathfrak{m}$ was obtained by A.A.~Vla\-di\-mirov and I.A.~Sheipak~\cite{VlaShe2010}
in the case ${\cal L}u\equiv-u''$. In~\cite{NazShei2012}, \cite{VlaShe2013-discrete} this result was generalized for the operators of arbitrary order.\footnote{In fact, in~\cite{VlaShe2010}, 
 \cite{VlaShe2013-discrete} sign-changing weights were also considered. Moreover, much more precise information on the spectrum structure was obtained in these papers. However, these results are not sufficient to receive exact small ball asymptotics for corresponding processes.} 
Using these results and Lemma~\ref{lm:log-P-slow-var-spect-asymp}, the logarithmic small ball asymptotics in $L_2[0,1;d\mathfrak{m}]$ was obtained in~\cite{NazShei2012} for several classical Green Gaussian processes.

\begin{theorem}[{\cite[Proposition 4.3]{NazShei2012}}] 
Let $\mathfrak{m}$ be a degenerate self-similar measure described above. Let $X$ be a Green Gaussian process, corresponding to the order~$2\ell$ ODO. Then
\begin{equation*}
\log {\mathbb P}\{\|X\|_{2,[0,1];\,d\mathfrak{m}}\leq\ep\}\sim
-\,\frac {n}{\log(\mathfrak{c}^{-1})}\,\log^2(\ep^{-1}),
\quad \ep\to 0,
\end{equation*}
where $\mathfrak{c}=\mathfrak{r}\cdot |I|^{2\ell-1}$.
\end{theorem}

Tensor products of Gaussian processes with eigenvalues counting functions, which are slowly varying at zero, were considered by A.I.~Karol' and A.I.~Nazarov~\cite{KarNaz2014}.\footnote{The small ball behavior of tensor product of some processes ${\cal R}_{C,2}$, see
\eqref{eq:spect-meas-smooth}, was investigated by T.~K\"uhn~\cite{Kuhn2011} in various $L_p$-norms. However, only two-sided estimates for small ball probabilities were obtained in~\cite{Kuhn2011}.}

\begin{theorem}[{\cite[Theorem 3.2]{KarNaz2014}}]
\label{thm:log-P-SVF}
Let ${\cal K}_1$ and ${\cal K}_2$ be positive (compact, self-adjoint) operators. Suppose that
\begin{align}
    {\cal N}_{{\cal K}_1}(\tau)\sim \Phi_1(\tau^{-1}), \qquad {\cal N}_{{\cal K}_2}(\tau)\sim \Phi_2(\tau^{-1}), \quad \tau\to0,
\end{align}
where $\Phi_1$ and $\Phi_2$ are unbounded and non-decreasing SVFs. Then we have
\begin{align}
{\cal N}_{\otimes}(\tau) \sim (\Phi_1\star\Phi_2)(\tau^{-1}),\qquad
\tau\to0, 
\end{align}
where 
\begin{align}
(\Phi\star\Psi)(s)
=
\int\limits_1^{s}\Phi(s\sigma^{-1})d\Psi(\sigma)\end{align}
is the so-called \textit{\textbf{asymptotic convolution}} of SVFs.\footnote{The properties of asymptotic convolution were investigated in~\cite[Section~2]{KarNaz2014}. In particular, the asymptotic convolution is connected with the Mellin convolution (see \eqref{eq:Mellin-conv}) by the relation 
$$
(\Phi\star\Psi)(s)=(\Phi*\Psi_1)(s),\qquad \mbox{where}\qquad \Psi_1(s)=s\Psi'(s).
$$
}
\end{theorem}

Using general results of Theorem~\ref{thm:log-P-SVF} and Lemma~\ref{lm:log-P-slow-var-spect-asymp}, the logarithmic $L_2$-small ball asymptotics were obtained in~\cite{KarNaz2014} for several fields-products with marginal processes described above.

\addcontentsline{toc}{section}{Acknowledgement}
\section*{Acknowledgement} 
We are deeply grateful to Prof. M.A. Lifshits who suggested us the idea of this survey, and to Prof. I.A. Sheipak who provided us with some important references. 

The work was supported by Russian Scientific Foundation, Grant 21-11-00047. The most part of this paper was written while Yu.P. was a postdoctoral fellow at IMPA. She thanks IMPA for creating excellent working conditions.

\appendix
\addcontentsline{toc}{section}{Appendix. $L_2$-small ball asymptotics for concrete processes}
\section*{Appendix. \newline $L_2$-small ball asymptotics for concrete processes}

\begin{tabular}{|m{5.5cm}|>{\centering\arraybackslash}m{2.2cm}|>{\centering\arraybackslash}m{2.2cm}|>{\centering\arraybackslash}m{2.2cm}|}
\hline 
Process 
&
Logarithmic asymptotics
&  
Asympt. up to constant
& 
Exact asymptotics
\\
\hline
\hline 
Wiener process, $W$
&&& \cite{CamMar1944}\phantom{8}
\\
\hline 
Brownian bridge, $B$ &&&\cite{AndDar1952}\phantom{88}
\\
\hline 
Ornstein-Uhlenbeck, $U$
&\cite{Li1992}&\cite{NazNik2004} &\cite{Naz2003-sharp};~~\cite{GaoHanLeeTor2003-Hadamard}\phantom{;}\phantom{;}\phantom{;}
\\
\hline 
Centered Wiener process, $\overline{W}$
&&&\cite{BegNikOrs2005}\phantom{8}
\\
\hline 
Centered Brownian bridge, $\overline{B}$
&&&\cite{BegNikOrs2005}\phantom{8}
\\
\hline 
Other Green Gaussian processes of the second order
&&&\cite{NikOrs2004};~~\cite{Naz2009-non-separated}; \cite{Pus2010-bogoliubov};~~\cite{Jin2014-phd}\phantom{8;}
\\
\hline 
Weighted Green Gaussian processes of the second order
&&\cite{Li1992};\phantom{\footnotemark[67]}~\cite{DehMar2003}\phantom{\footnotemark[67];}&
\cite{Naz2003-sharp};~~\cite{GaoHanLeeTor2003-Hadamard};\phantom{8} \cite{NikKha2006};~~\cite{NazPus2009};
\cite{Pus2010-bogoliubov};~~\cite{BarIgl2011};\phantom{88}
\cite{NikPus2013};~~\cite{NazNik2018-mixed}\phantom{;}
\\
\hline
\hline
$m$-times integrated  Wiener process,  $W_m^{[\beta_1,\ldots,\beta_m]}$
&
\cite{KhoShi1998}\footnotemark[67];~\cite{ChenLi2003}\phantom{;}
&\cite{BegNikOrs2005}\footnotemark[67];~\cite{GaoHanTor2003-integrated}\footnotemark[68]; \cite{NazNik2004}
&\cite{GaoHanLeeTor2003-Hadamard};\phantom{8}~~\cite{GaoHanLeeTor2004};\phantom{8} \cite{Naz2003-sharp} 
\\
\hline
$m$-times integrated Bro\-wnian bridge, $B_m^{[\beta_1,\ldots,\beta_m]}$ 
& 
&
\cite{BegNikOrs2005}\footnotemark[67];~\cite{NazNik2004}\phantom{;}
& 
\cite{GaoHanLeeTor2003-Hadamard};\phantom{8}~~\cite{Naz2003-sharp}\phantom{;}
\\
\hline 
$m$-times integrated Ornstein--Uhlenbeck, $U_m^{[\beta_1,\ldots,\beta_m]}$ 
&
&
\cite{NazNik2004}
&
\cite{Naz2003-sharp}
\\
\hline
Other Green Gaussian processes of the higher order
&&\cite{NazNik2004};\phantom{;}~\cite{BegNikOrs2005}\phantom{8;} &
\cite{Naz2003-sharp};~~\cite{Pus2008-Matern}; 
\cite{Naz2009-non-separated};~~\cite{Pus2010-bogoliubov}; 
\cite{Pus2011-Matern};~~\cite{NazNik2018-mixed}\phantom{;}
\\
\hline
Weighted Green Gaussian processes of the higher order
&&&\cite{NazPus2009};~~\cite{NazPus2014-comparison}\phantom{;}
\\
\hline
\end{tabular}

\footnotetext[67]{The case $m=1$.}
\footnotetext[68]{The exponent of the power term was determined up to a constant, see~\eqref{eq:small-ball-integrated-w-gao}.}

\newpage

\begin{tabular}{|m{5.5cm}|>{\centering\arraybackslash}m{2.2cm}|>{\centering\arraybackslash}m{2.2cm}|>{\centering\arraybackslash}m{2.2cm}|}
\hline 
Process 
&
Logarithmic asymptotics
&  
Asympt. up to constant
& 
Exact asymptotics
\\
\hline 
\hline 
Detrended Green processes
&&\cite{AiLiLiu2012}\phantom{88}
&
\cite{KirNik2014};\phantom{8}~~\cite{Pet2017}\phantom{;}
\\
\hline
Finite-dimensional perturbations\footnotemark[69]
&&&\cite{Naz2010-onedimpert};~~\cite{NazPet2015-SVFpert}; \cite{Pet2020}
\\
\hline
\hline
Fractional Brownian motion (fBm), $W^H$
&\cite{Bro2003}\phantom{8}
&\cite{ChiKlep2018}\phantom{8}
&
\\
\hline 
$m$-times integrated fBm
&
\cite{NazNik2005-frac}
&
\cite{ChiKlepMar2018}\footnotemark[70];~\cite{ChiKlepNazRast2023}\phantom{8;}
&
\\
\hline 
Other fractional-Green processes
&
&
\cite{ChiKlepMar2020};~~~\cite{Naz2021-frac}; \cite{ChiKlepNazRast2023}\phantom{8}
&
\\
\hline
Riemann--Liouville process and bridge
&
&
\cite{ChiKlep2021}\phantom{8}
&
\\
\hline 
Other fractional processes and fields
&\cite{NazNik2005-frac};~~\cite{Pus2011-Matern}; \cite{HongNazLif2016};\phantom{8}~~\cite{LifNaz2018}\phantom{;}&&
\\
\hline
Smooth processes
& \cite{GaoLi2007-log};\phantom{8}~~\cite{Naz2009-log}; \cite{AurGaoKuhnLiShao2013}\phantom{88}
&&
\\
\hline
Fractional processes with variable Hurst index
&\cite{KarNaz2021-Multifrac};\phantom{8}~~\cite{Kar2022}\phantom{8;}&&
\\
\hline
Processes with self-similar measure
&\cite{Naz2004-logselfsim};~~~\cite{NazShei2012}; 
\cite{FreRast2020-sib}\footnotemark[71];~\cite{Rast2020-selfsimilar}\phantom{;}
&&
\\
\hline 
\hline 
Brownian sheet
&
\cite{Csa1982}\phantom{8}
&
&
\cite{FillTor2004}\phantom{8}
\\
\hline
Integrated Brownian sheet
&\cite{KarNazNik2008-Tensor}\phantom{8}
&
&
\cite{FillTor2004}\phantom{8}
\\
\hline
Other fields-products
&\cite{KarNazNik2008-Tensor};\phantom{8}~~\cite{GaoLi2007-log};\phantom{8} 
\cite{GaoLi2007-slepian};\phantom{8}~~\cite{KarNaz2014};\phantom{8}
\cite{Kuhn2011}\footnotemark[71];~\cite{Rast2018-tensor}\phantom{;}
&
&
\cite{Roz2019-fields};~~\cite{Roz2019}\phantom{;}
\\
\hline

\end{tabular}
\footnotetext[69]{Perturbations from the kernel of original Gaussian random vector.}
\footnotetext[70]{The case $m=1$.}
\footnotetext[71]{Two-sided estimates.}

\newpage
\addcontentsline{toc}{section}{References}
\small

\end{document}